\documentclass[11pt, leqno]{article}

\usepackage{amssymb,amsmath,amsthm,mathrsfs}
\usepackage{graphicx}
\usepackage{times}

\leftmargini=\baselineskip

\newtheoremstyle{localthm}
	{5pt} 
	{5pt} 
	{\sl} 
	{} 
	{\bf} 
	{{\rm.}} 
	{.7em} 
	{} 

\theoremstyle{localthm}
\newtheorem{Theorem}{Theorem}[section]
\newtheorem{Corollary}[Theorem]{Corollary}
\newtheorem{Proposition}[Theorem]{Proposition}
\newtheorem{Lemma}[Theorem]{Lemma}

\newtheoremstyle{localrem}
	{5pt} 
	{5pt} 
	{\rm} 
	{} 
	{\bf} 
	{{\rm.}} 
	{.7em} 
	{} 

\theoremstyle{localrem}
\newtheorem{Remark}[Theorem]{Remark}

\newcommand{\R}{\mathbb{R}}
\renewcommand{\d}{\mathrm{d}}
\newcommand{\veps}{\varepsilon}
\newcommand{\Ex}{\operatorname{\mathbb{E}}}

\begin{document}
\addtolength{\baselineskip}{0.2\baselineskip}

\title{Various New Inequalities for Beta Distributions}
\author{Alexander Henzi and Lutz D{\"u}mbgen\\
	ETH Z{\"u}rich and University of Bern}
\date{August 2023}

\maketitle

\begin{abstract}
This note provides some new inequalities and approximations for beta distributions, including tail inequalities, exponential inequalities of Hoeffding and Bernstein type, Gaussian inequalities and approximations.

\bigskip

Parts of Version~8 of this manuscript have been published as ``Some New Inequalities for Beta Distributions'' in \textsl{Statistics \& Probability Letters \textbf{195}} (2023). The present Version~9 provides even stronger inequalities in Sections~\ref{sec:Segura} and \ref{sec:Exponential}, and some references have been updated.
\end{abstract}


\paragraph{MSC2020 subject classifications:} Primary 62E17; 60E05; 60E15, Secondary 33B20.

\section{Introduction}

Beta distributions play an important role in statistics and probability theory (Gupta and Nadarajah, 2004), and they occur in various scientific fields (Skorski, 2023). A frequent obstacle in problems involving beta distributions is the lack of analytic expressions for their distribution function, the normalized incomplete beta function. Therefore one often resorts to inequalities and approximations, as, for example, in the proofs of Dimitriadis et al.\ (2022, Theorem~4.1) and D{\"u}mbgen and Wellner (2023, Lemma~S.8).

This paper provides some new inequalities for the beta distribution $\mathrm{Beta}(a,b)$ with parameters $a, b > 0$, its distribution function $B_{a,b}$, survival function $\bar{B}_{a,b} = 1 - B_{a,b}$ and density function $\beta_{a,b}$ on $[0,1]$. The latter is given by
\[
	\beta_{a,b}(x) \ := \ B(a,b)^{-1} x^{a-1} (1 - x)^{b-1} , \quad x \in (0,1) ,
\]
wheere $B(a,b) := \int_0^1 x^{a-1}(1 - x)^{b-1} \, \d x = \Gamma(a) \Gamma(b)/\Gamma(a+b)$, and $\Gamma(\cdot)$ denotes the gamma function. In Section~\ref{sec:Segura}, we refine the lower and upper bounds for $B_{a,b}$ and $\bar{B}_{a,b}$ by Segura~(2016) which are particularly accurate in the tails of $\mathrm{Beta}(a,b)$. As a by-product we obtain refinements of Segura's~(2014) bounds for the gamma distribution and survival functions. In Section~\ref{sec:Exponential} we present new exponential inequalities which are stronger than previously known inequalities of D{\"u}mbgen~(1998), Marchal and Arbel~(2017) or Skorski~(2023) for a wide range of parameters $(a,b)$. Section~\ref{sec:Gaussian.tail.inequalities} presents inequalities for $B_{a,b}$ and $\bar{B}_{a,b}$ in terms of Gaussian distribution functions. Finally, Section~\ref{sec:Gaussian.approximation.Beta(a,a)} discusses the approximation of the symmetric distribution $\beta_{a,a}$ by Gaussian densities with mean $1/2$ in the spirit of D{\"u}mbgen et al.\ (2021). Most proofs are deferred to Section~\ref{sec:Proofs}.

\section{Sharp tail inequalities}
\label{sec:Segura}

In what follows, let $p := a/(a+b)$, the mean of $\mathrm{Beta}(a,b)$. In a general setting including noncentral beta distributions, Segura (2016, inequalities (27), (29), (30)) uses extensions of l'Hopital's rule to derive inequalities for $B_{a,b}$ and $\bar{B}_{a,b}$. For symmetry reasons, we only consider $B_{a,b}$, because $\bar{B}_{a,b}(x) = B_{b,a}(1-x)$. Elementary calculations reveal that
\[
	B_{a,b}(x) \
	= \ \frac{x^a}{a B(a,b)} \, Q_{a,b}(x)
\]
with
\[
	Q_{a,b}(x) \ := \ a \int_0^1 y^{a-1} (1 - xy)^{b-1} \, \d y .
\]
This auxiliary function $Q_{a,b}$ is well-defined and smooth on $(-\infty,1]$, and elementary calculations reveal that \[
	Q_{a,b}^{}(0) \ = \ 1 , \quad
	Q_{a,b}'(0) \ = \ - \frac{a(b-1)}{a+1} , \quad
	Q_{a,b}''(0) \ = \ \frac{a(b-1)(b-2)}{a+2} .
\]
More generally,
\[
	Q_{a,b}(x) \ = \ \sum_{k=0}^\infty \frac{a [b-1]_k (-x)^k}{(a+k) k!}
\]
with the descending Pochhammer symbols $[\cdot]_k$; in fact, $Q_{a,b}$ coincides with the hypergeometric function ${}_2F_1(1-b,a;a+1;\cdot)$, see Olver et al.\ (2010). Interesting special cases are:
\[
	Q_{a,1}(x) \ = \ 1 , \quad
	Q_{a,2}(x) \ = \ 1 - \frac{ax}{a+1} , \quad
	Q_{a,3}(x) \ = \ 1 - \frac{2ax}{a+1} + \frac{ax^2}{a+2} .
\]
Segura~(2016) showed that for $x \in [0,1)$,
\begin{equation}
\label{ineq:Segura.left}
	(1 - x)^b (1 + c_{a,b} x) \
	\le \ Q_{a,b}(x) \
	\le \ \frac{(1 - x)^b}{(1 - x/p)^+} , \\
\end{equation}
where $c_{a,b} := (a+b)/(a+1)$. Numerical examples reveal that these inequalities are rather accurate if $x \ll p$. Indeed, both bounds are equal to $1$ for $x = 0$, and the lower bound has derivative $-b + c_{a,b} = Q_{a,b}'(0)$ for $x = 0$. The upper bound is less accurate, because its derivative for $x = 0$ equals $- b + 1/p > Q_{a,b}'(0)$. Our first contribution is an improvement of Segura's bounds.

\begin{Theorem}
\label{thm:Segura1}
For $x \in [0,1)$, let
\begin{align*}
	Q_{a,b}^{[S,1]}(x) \
	&:= \ (1 - x)^b \Bigl( 1 + \frac{c_{a,b} x}{1 - x} \Bigr) , \\
	Q_{a,b}^{[S,2]}(x) \
	&:= \ \begin{cases}
		\displaystyle
		\frac{(1 - x)^b}{1 - c_{a,b}x} & \text{if} \ x \in [0,p] , \\
		(a+1)(1 - p)^b & \text{if} \ x \in [p,1) .
	\end{cases}
\end{align*}
Then,
\[
	Q_{a,b}^{}(x) \ \begin{cases}
		\le \ Q_{a,b}^{[S,1]}(x) & \text{if} \ b \le 1 , \\[0.5ex]
		\ge \ Q_{a,b}^{[S,1]}(x) & \text{if} \ b \ge 1 , \\[1ex]
		\ge \ Q_{a,b}^{[S,2]}(x) & \text{if} \ b \le 1 , \\[0.5ex]
		\le \ Q_{a,b}^{[S,2]}(x) & \text{if} \ b \ge 1 .
	\end{cases}
\]
\end{Theorem}

\begin{Remark}
\label{rem:Segura1}
Note that $Q_{a,1}^{[S,1]} \equiv Q_{a,1}^{[S,2]} \equiv Q_{a,1}^{} \equiv 1$. In general, for $\ell= 1,2$, the function $Q(x) := Q_{a,b}^{[S,\ell]}(x)$ satisfies $Q(0) = 1$ and $Q'(0) = -b + c_{a,b} = Q_{a,b}'(0)$.
\end{Remark}

\begin{figure}
\centering
\includegraphics[width=0.85\textwidth]{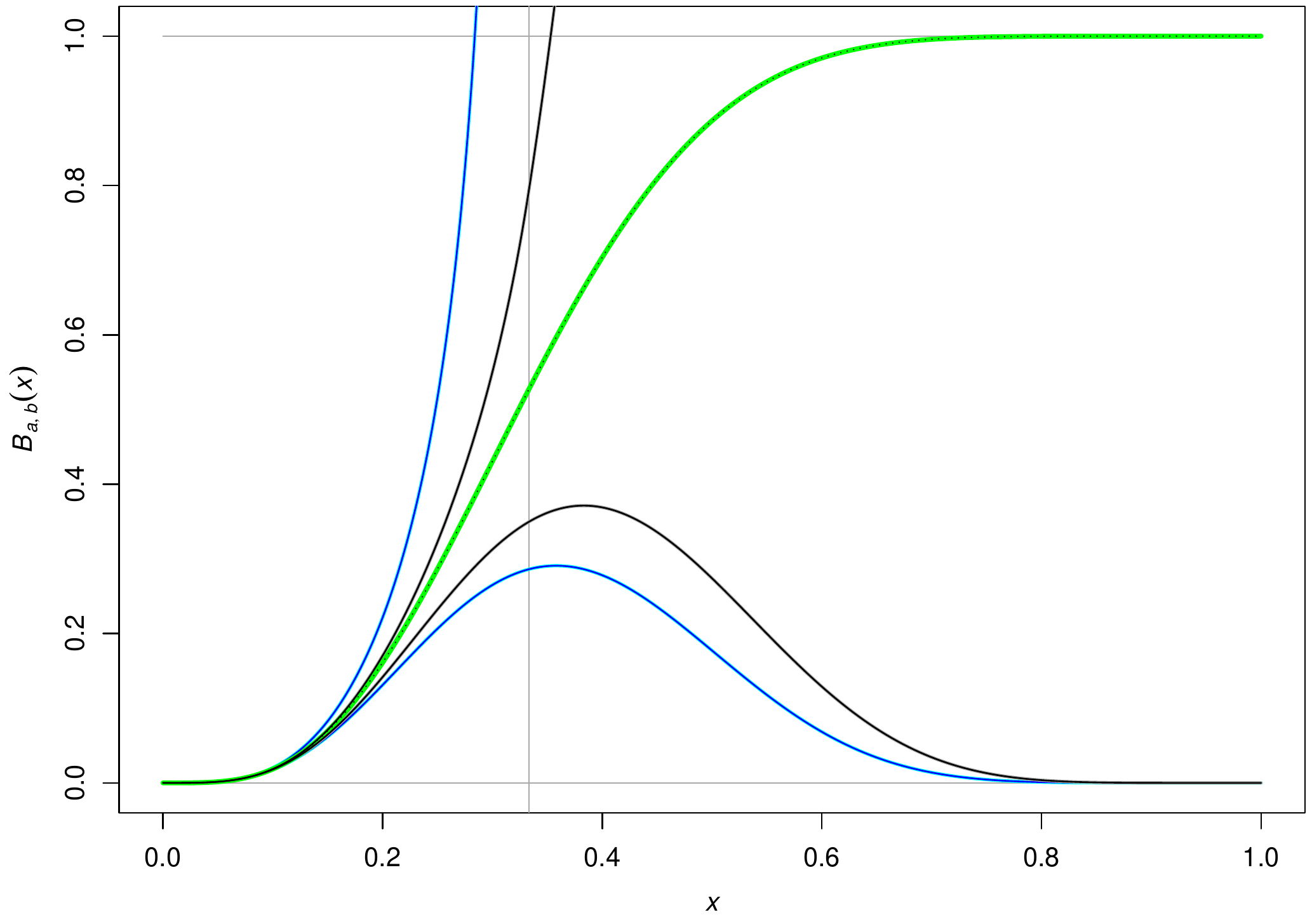}

\includegraphics[width=0.85\textwidth]{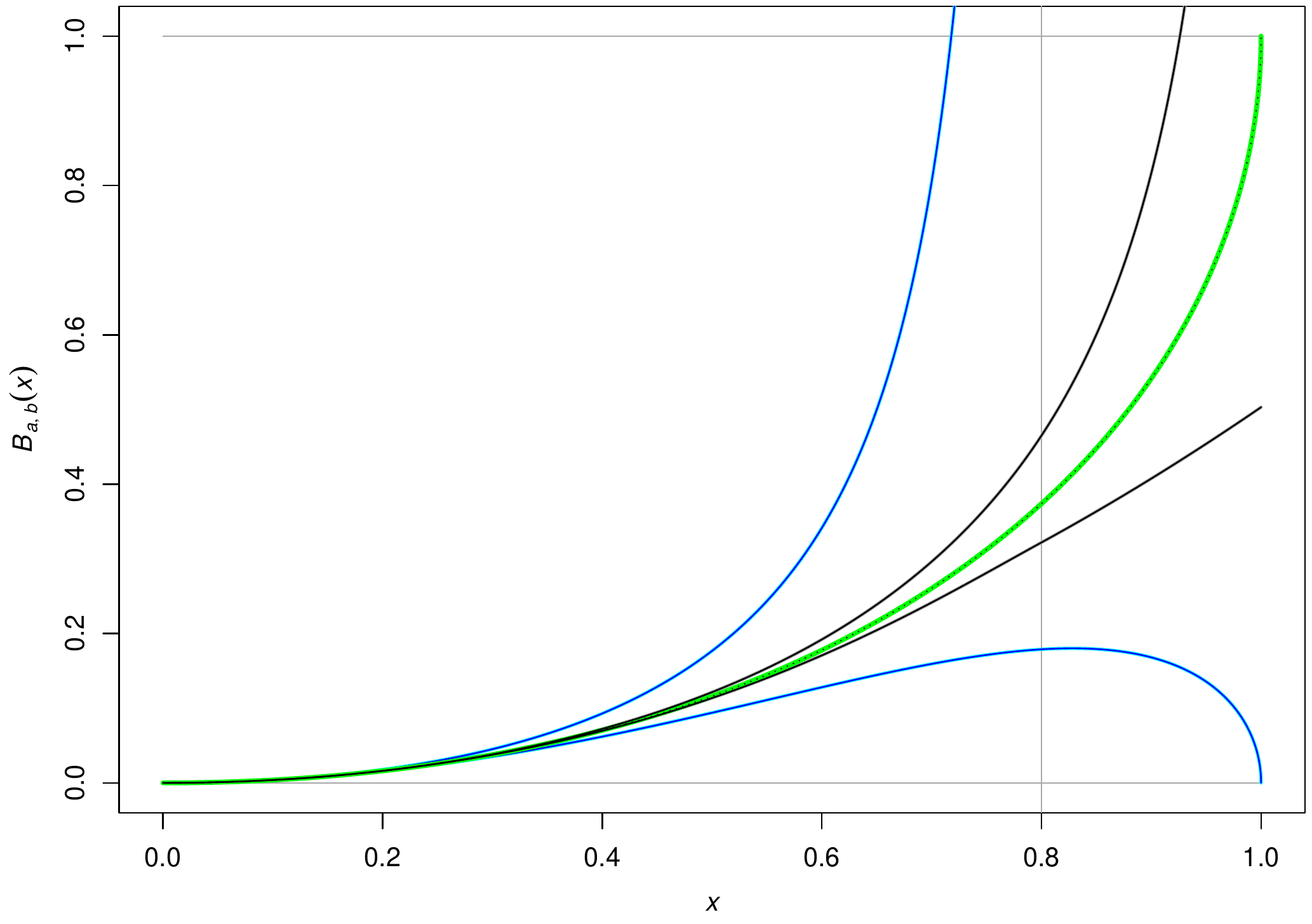}

\caption{Inequalities for $B_{a,b}$ when $(a,b) = (4,8)$ (upper panel) and $(a,b) = (2,0.5)$ (lower panel). The green line shows $B_{a,b}$, the blue lines are Segura's (2016) bounds resulting from \eqref{ineq:Segura.left}, and the black lines are the bounds via Theorem~\ref{thm:Segura1}. The vertical line indicates the mean $p$.}
\label{fig:Segura}
\end{figure}

Figure~\ref{fig:Segura} illustrates the bounds for $B_{a,b}$ resulting from \eqref{ineq:Segura.left} and Theorem~\ref{thm:Segura1} in case of $(a,b) = (4,8), (2,0.5)$.

The next result provides alternative bounds for $Q_{a,b}$.

\begin{Theorem}
\label{thm:Segura2}
For $x \in (0,1)$ let
\begin{align*}
	Q_{a,b}^{[1]}(x) \
	&:= \ \Bigl( 1 - \frac{ax}{a+1} \Bigr)^{b-1} , \\
	Q_{a,b}^{[2]}(x) \
	&:= \ (1 - x)^{b-1} + \frac{(b-1)x}{a+1} \Bigl( 1 - \frac{(a+1)x}{a+2} \Bigr)^{b-2} , \\
	Q_{a,b}^{[3]}(x) \
	&:= \ \frac{1}{(a+1)^2} + \frac{a(a+2)}{(a+1)^2} \Bigl( 1 - \frac{(a+1)x}{a+2} \Bigr)^{b-1} , \\
	Q_{a,b}^{[4]}(x) \
	&:= \ \frac{a(1 - x)^{b-1} + 1}{a+1}
		- \frac{a(b-1)(b-2)x^2}{2(a+1)(a+2)} \Bigl( 1 - \frac{2(a+2)x}{3(a+3)} \Bigr)^{b-3} , \\
	Q_{a,b}^{[5]}(x) \
	&:= \ \frac{a(1 - x)^{b-1} + 1}{a+1}
		- \frac{a(b-1)(b-2)x^2 \bigl( a + 5 + 2(a+2)(1-x)^{b-3} \bigr)}{6(a+1)(a+2)(a+3)} .
\end{align*}
Then,
\[
	Q_{a,b}^{}(x) \ \begin{cases}
		\ge \ Q_{a,b}^{[1]}(x) & \text{if} \ b \in (0,1] \cup [2,\infty) , \\[0.5ex]
		\le \ Q_{a,b}^{[1]}(x) & \text{if} \ b \in [1,2] , \\[1ex]
		\le \ Q_{a,b}^{[2]}(x) & \text{if} \ b \in (0,1] \cup [2,3] , \\[0.5ex]
		\ge \ Q_{a,b}^{[2]}(x) & \text{if} \ b \in [1,2] \cup [3,\infty) , \\[1ex]
		\ge \ Q_{a,b}^{[3]}(x) & \text{if} \ b \in (0,1] \cup [2,3] , \\[0.5ex]
		\le \ Q_{a,b}^{[3]}(x) & \text{if} \ b \in [1,2] \cup [3,\infty) , \\[1ex]
		\le \ Q_{a,b}^{[4]}(x) & \text{if} \ b \in (0,1] \cup [2,3] \cup [4,\infty) , \\[0.5ex]
		\ge \ Q_{a,b}^{[4]}(x) & \text{if} \ b \in [1,2] \cup [3,4] , \\[1ex]
		\ge \ Q_{a,b}^{[5]}(x) & \text{if} \ b \in (0,1] \cup [2,3] \cup [4,\infty)  , \\[0.5ex]	
		\le \ Q_{a,b}^{[5]}(x) & \text{if} \ b \in [1,2] \cup [3,4] .
	\end{cases} 
\]
\end{Theorem}

\begin{Remark}
\label{rem:Segura2}
Note that the inequalities in terms of $Q_{a,b}^{[\ell]}$ are equalities for
\[
	\begin{cases}
		\ell = 1 \ \text{and} \ b = 1,2 , \\
		\ell = 2,3 \ \text{and} \ b = 1,2,3 , \\
		\ell = 4,5 \ \text{and} \ b = 1,2,3,4 .
	\end{cases}
\]
Note further that for $\ell = 1,\ldots,5$, the function $Q := Q_{a,b}^{[\ell]}$ satisfies $Q(0) = 1$, $Q'(0) = Q_{a,b}'(0)$, while
\[
	Q''(0) \ = \ \begin{cases}
		\displaystyle
		\frac{a(a+2)}{(a+1)^2} Q_{a,b}''(0) & \text{for} \ \ell = 1 , \\[1.5ex]
		Q_{a,b}''(0) & \text{for} \ \ell = 2,3,4,5 .
	\end{cases}
\]
\end{Remark}

\begin{Remark}
\label{rem:Segura12}
The lower bounds for $Q_{a,b}$ resulting from Theorem~\ref{thm:Segura2} are stronger than the ones from Theorem~\ref{thm:Segura1}. Precisely, it is shown in Section~\ref{sec:Proofs} that the following inequalities hold true for $x \in (0,1]$:
\begin{align*}
	Q_{a,b}^{[S,2]}(x) \ &< \ Q_{a,b}^{[1]}(x) \quad \text{if} \ b \in (0,1) , \\
	Q_{a,b}^{[S,1]}(x) \ &< \ Q_{a,b}^{[2]}(x) \quad \text{if} \ b \in (1,2] , \\
	Q_{a,b}^{[S,1]}(x) \ &< \ Q_{a,b}^{[1]}(x) \quad \text{if} \ b \in [2,\infty) .
\end{align*}
Moreover, the upper bounds for $Q_{a,b}$ resulting from Theorem~\ref{thm:Segura2} are stronger than the ones from Theorem~\ref{thm:Segura1} if $b \in (0,3]$. Precisely, the following inequalities are derived in Section~\ref{sec:Proofs} for $x \in (0,1]$:
\begin{align*}
	Q_{a,b}^{[S,1]}(x) \ &> \ Q_{a,b}^{[2]}(x) \quad \text{if} \ b \in (0,1) , \\
	Q_{a,b}^{[S,2]}(x) \ &> \ Q_{a,b}^{[1]}(x) \quad \text{if} \ b \in (1,2] , \\
	Q_{a,b}^{[S,2]}(x) \ &> \ Q_{a,b}^{[2]}(x) \quad \text{if} \ b \in [2,3] \ \text{and} \ x \le p .
\end{align*}
\end{Remark}

\begin{figure}
\centering
\includegraphics[width=0.85\textwidth]{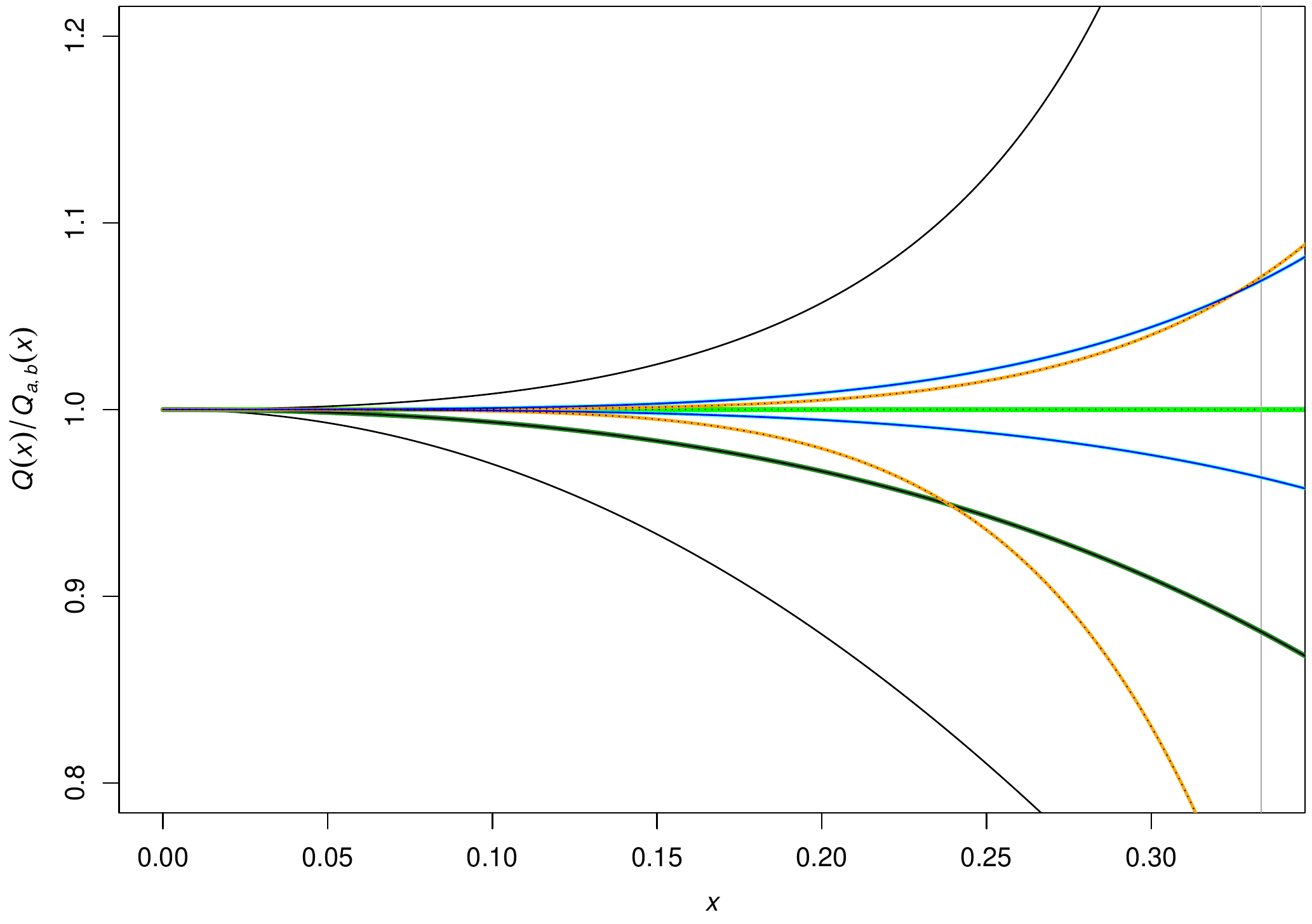}

\includegraphics[width=0.85\textwidth]{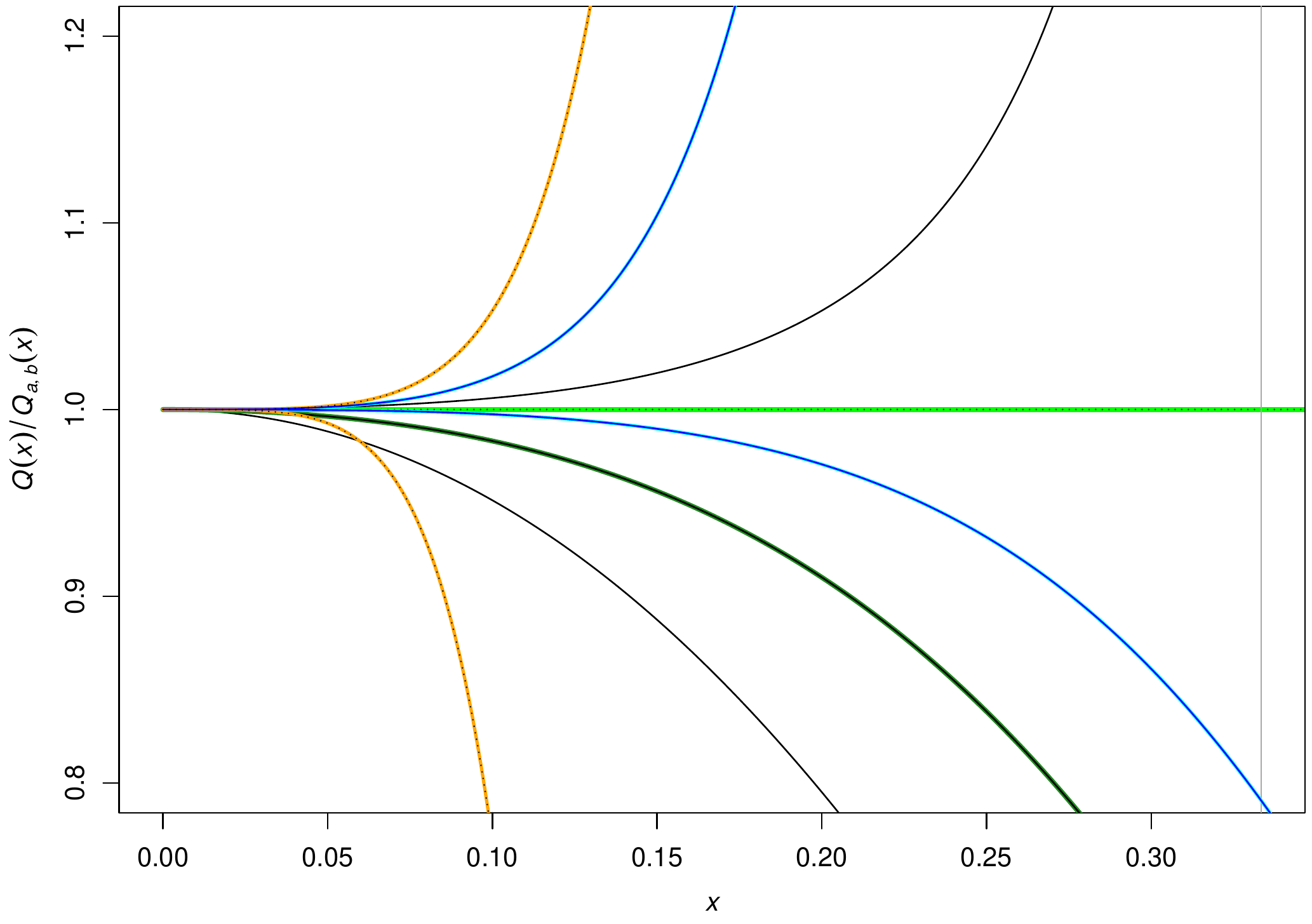}

\caption{Inequalities for $Q_{a,b}$ when $(a,b) = (4,8)$ (upper panel) and $(a,b) = (12,24)$ (lower panel) on $[0,p]$. We depict the ratios $Q_{a,b}^{[S,\ell]}/Q_{a,b}^{}$ for $\ell = 1,2$ (thin black) and $Q_{a,b}^{[\ell]}/Q_{a,b}^{}$ for $\ell = 1$ (dark green), $\ell = 2,3$ (blue) and $\ell = 4,5$ (orange, black dashed).}
\label{fig:Segura2A}
\end{figure}

\begin{figure}
\centering
\includegraphics[width=0.85\textwidth]{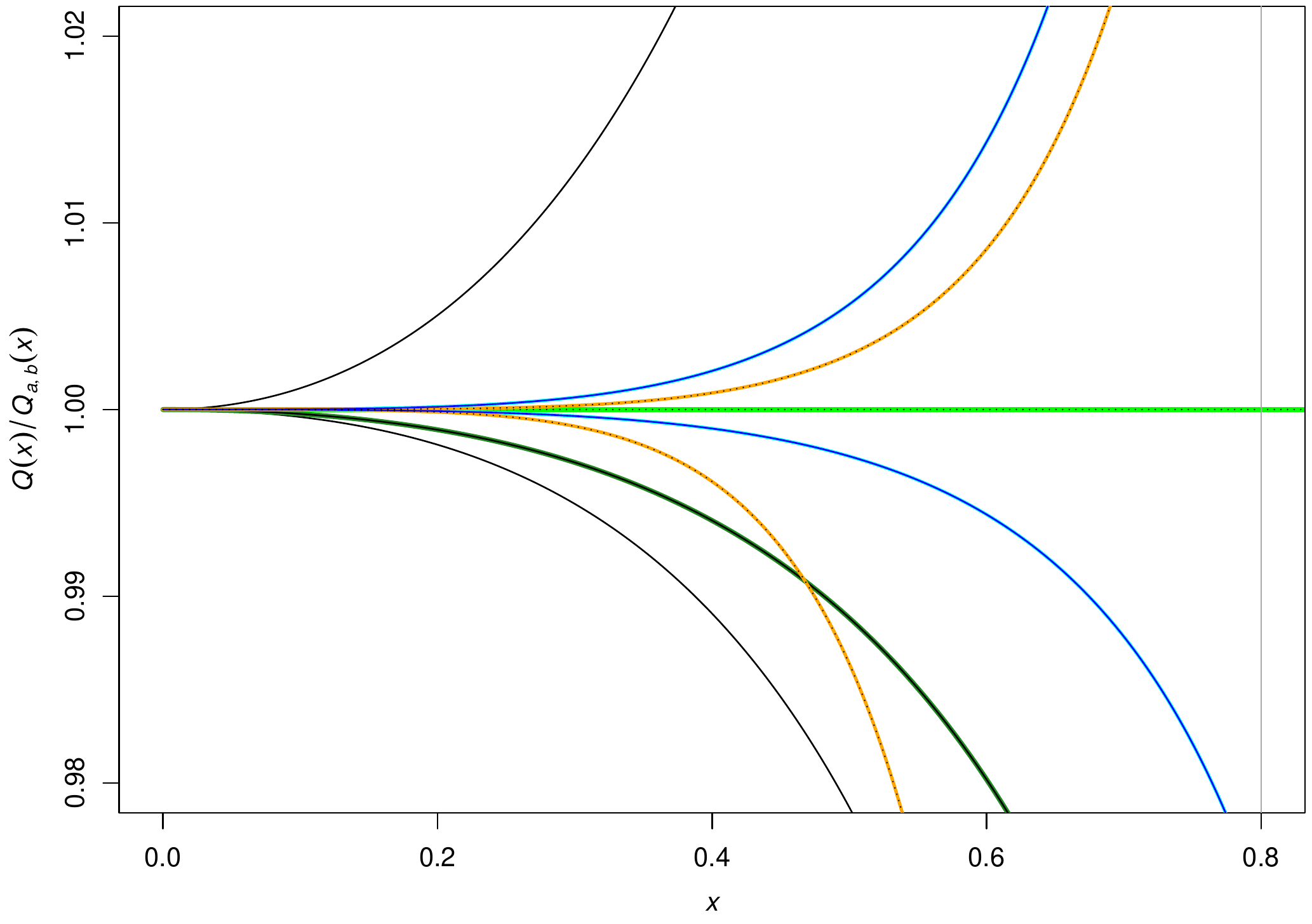}

\includegraphics[width=0.85\textwidth]{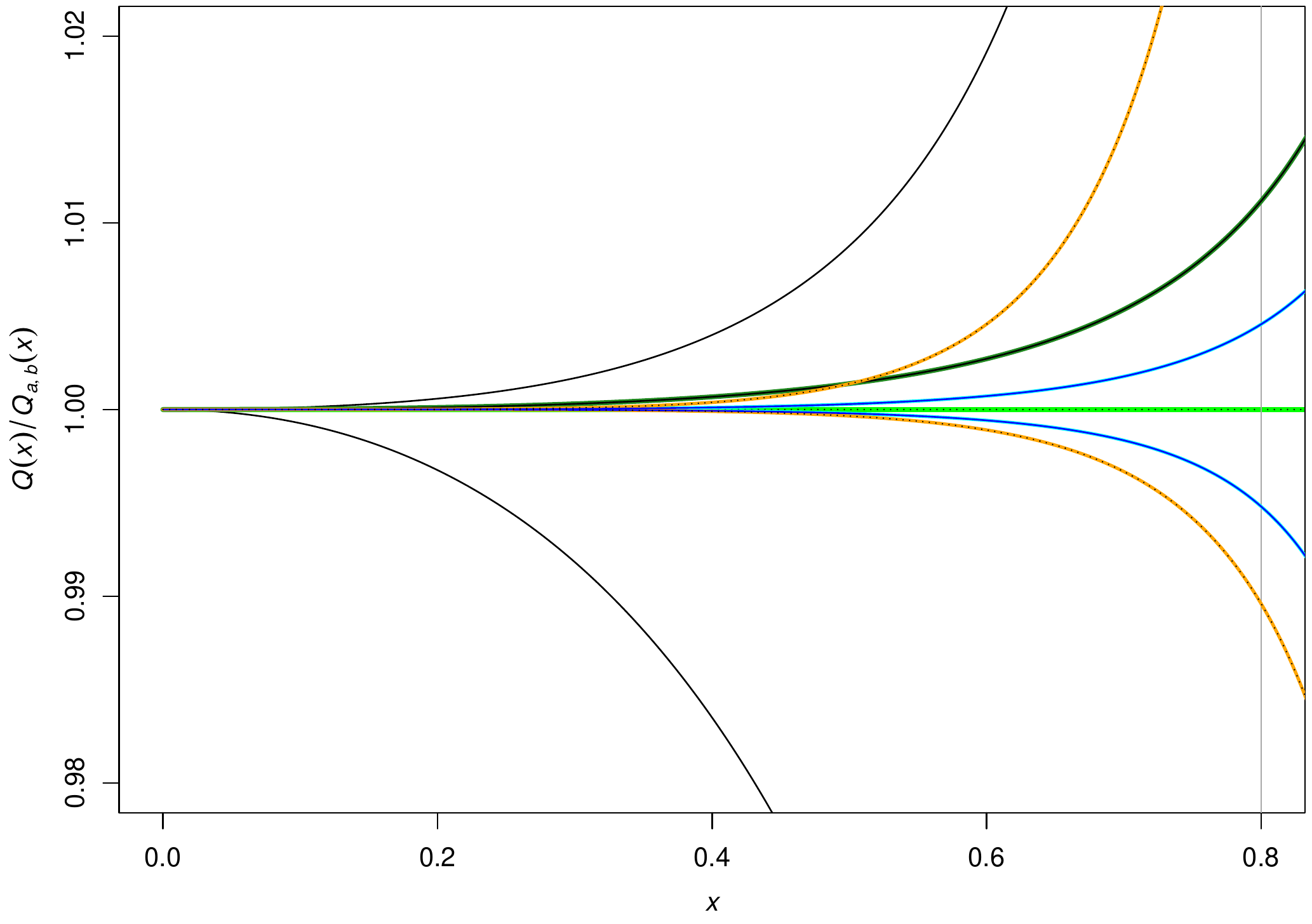}

\caption{Inequalities for $Q_{a,b}$ when $(a,b) = (2,0.5)$ (upper panel) and $(a,b) = (6,1.5)$ (lower panel) on $[0,p]$. We depict the ratios $Q_{a,b}^{[S,\ell]}/Q_{a,b}^{}$ for $\ell = 1,2$ (thin black) and $Q_{a,b}^{[\ell]}/Q_{a,b}^{}$ for $\ell = 1$ (dark green), $\ell = 2,3$ (blue) and $\ell = 4,5$ (orange, black dashed).}
\label{fig:Segura2B}
\end{figure}

Figures~\ref{fig:Segura2A} and \ref{fig:Segura2B} illustrate the bounds for $Q_{a,b}$ resulting from Theorems~\ref{thm:Segura1} and \ref{thm:Segura2} for various pairs $(a,b)$. In all cases one sees the ratio $Q/Q_{a,b}$ on $[0,p]$, where $Q$ stands for one of seven bounds.

\paragraph{Gamma distributions.}
There is a rich literature about inequalities for gamma distribution and survival functions, see, for instance, Qi and Mei~(1999), Neuman~(2013), Segura~(2014) and Pinelis~(2020). We just illustrate that our bounds in Theorems~\ref{thm:Segura1} and \ref{thm:Segura2} yield a connection to that literature. It is well-known that for a random variable $X_{a,b} \sim \mathrm{Beta}(a,b)$, the rescaled variable $b X_{a,b}$ converges in distribution to a gamma random variable with shape parameter $a$ and scale parameter $1$ as $b \to \infty$. Denoting the corresponding distribution and survival function with $G_a$ and $\bar{G}_a = 1 - G_a$, respectively, we have $G_a(x) = \lim_{b\to \infty} B_{a,b}(x/b)$, and one can deduce from Theorems~\ref{thm:Segura1} and \ref{thm:Segura2} the following bounds.

\begin{Corollary}
\label{cor:Segura}
For $a, x > 0$, define
\begin{align*}
	Q_a^{[S]}(x) \
	&:= \ \frac{(a+1)}{(a+1-x)^+} \, e_{}^{-x} , \\
	Q_a^{[1]}(x) \
	&:= \ e_{}^{-ax/(a+1)} , \\
	Q_a^{[2]}(x) \
	&:= \ e_{}^{-x} + \frac{x}{a+1} \, e_{}^{-(a+1)x/(a+2)} , \\
	Q_a^{[3]}(x) \
	&:= \ \frac{1}{(a+1)^2} + \frac{a(a+2)}{(a+1)^2} \, e_{}^{-(a+1)x/(a+2)} , \\
	Q_a^{[4]}(x) \
	&:= \ \frac{ae^{-x} + 1}{a+1} - \frac{ax^2}{2(a+1)(a+2)} \, e_{}^{-2(a+2) x/[3(a+3)]} , \\
	Q_a^{[5]}(x) \
	&:= \ \frac{ae^{-x} + 1}{a+1} - \frac{ax^2 [a+5 + 2(a+2) e^{-x}]}{6(a+1)(a+2)(a+3)} .
\end{align*}
Then
\[
	\frac{a \Gamma(a)}{x^a} \, G_a^{}(x) \ \begin{cases}
		\le \ \min \bigl\{ Q_a^{[S]}(x), Q_a^{[3]}(x), Q_a^{[4]}(x) \bigr\} , \\
		\ge \ \max \bigl\{ Q_a^{[1]}(x), Q_a^{[2]}(x), Q_a^{[5]}(x) \bigr\} .
	\end{cases}
\]
Moreover, define
\begin{align*}
	\bar{Q}_a^{[S]}(x) \
	&:= \ 1_{[x \le a]}^{} a^a + 1_{[x >a]} \frac{x^a}{x-a+1} , \\
	\bar{Q}_a^{[1]}(x) \
	&:= \ (x + 1)^{a-1} , \\
	\bar{Q}_a^{[2]}(x) \
	&:= \ x^{a-1} + (a-1) (x + 1)^{a-2} .
\end{align*}
Then
\[
	\Gamma(a) e_{}^x \, \bar{G}_a^{}(x) \ \begin{cases}
		\ge \ \bar{Q}_a^{[1]}(x) & \text{if} \ a \in (0,1] \cup [2,3] , \\
		\le \ \bar{Q}_a^{[1]}(x) & \text{if} \ a \in [1,2] , \\[1ex]
		\le \ \bar{Q}_a^{[2]}(x) & \text{if} \ a \in (0,1] \cup [2,3] , \\
		\ge \ \bar{Q}_a^{[2]}(x) & \text{if} \ a \in [1,2] \cup [3,\infty) , \\[1ex]
		\le \ \bar{Q}_a^{[S]}(x) & \text{if} \ a \in [3,\infty) .
	\end{cases}
\]
\end{Corollary}

Some of these bounds are known already or refinements of results in the aforementioned literature, notably Neuman~(2013, Theorem~4.1) and Segura~(2014, Theorem~10). Note also that for integers $a \ge 1$,
\[
	\Gamma(a) e_{}^x \, \bar{G}_a^{}(x) \ = \ \Gamma(a) \sum_{i=0}^{a-1} \frac{x^i}{i!} .
\]
The bounds in Corollary~\ref{cor:Segura} reproduce this equality for $a = 1,2,3$.

\section{Exponential inequalities}
\label{sec:Exponential}

Although the upper bounds in Theorems~\ref{thm:Segura1} and \ref{thm:Segura2} are numerically rather accurate in the tails, they can diverge to $\infty$ at $x = p$ as $a,b \to \infty$. Moreover, it is sometimes desirable to have bounds for $\log B_{a,b}(x)$ and $\log \bar{B}_{a,b}(x)$ in terms of simple, maybe rational, functions of $x$. Numerous exponential tail inequalities for $B_{a,b}$ and $\bar{B}_{a,b}$ have been derived already. We start with one particular result of D{\"u}mbgen (1998, Proposition~2.1). For $x \in [0,1]$ let
\[
	K(p,x) \
	:= \ p \log \Bigl( \frac{p}{x} \Bigr)
		+ (1 - p) \log \Bigl( \frac{1 - p}{1 - x} \Bigr)
		\in [0,\infty] .
\]
This function $K(p,\cdot)$ is strictly convex with minimum $K(p,p) = 0$. For arbitrary $x \in [0,1]$,
\begin{equation}
\label{ineq:Duembgen}
	\Bigl( \frac{x}{p} \Bigr)^a \Bigl( \frac{1 - x}{1 - p} \Bigr)^b
		= \exp \bigl( - (a+b) K(p,x) \bigr) \
	\ge \ \begin{cases}
		B_{a,b}(x)
			& \text{if} \ x \le p , \\
		\bar{B}_{a,b}(x)
			& \text{if} \ x \ge p .
	\end{cases}
\end{equation}
These inequalities can be improved substantially. The starting point is a rather general inequality.

\begin{Proposition}
\label{prop:Beta}
Let $q \in (0,1)$. For $x \in [0,q]$,
\[
	B_{a,b}(x) \ \le \ B_{a,b}(q)
		\Bigl( \frac{x}{q} \Bigr)^a \Bigl( \frac{1-x}{1-q} \Bigr)^{c_\ell}
	\quad\text{with} \ \ c_\ell := \frac{a (b-1)^+}{a+1} ,
\]
and for $x \in [q,1]$,
\[
	\bar{B}_{a,b}(x) \ \le \ \bar{B}_{a,b}(q)
		\Bigl( \frac{x}{q} \Bigr)^{c_r} \Bigl( \frac{1-x}{1-q} \Bigr)^b
	\quad\text{with} \ \ c_r := \frac{b (a-1)^+}{b+1} .
\]
\end{Proposition}

Combining this proposition with well-known inequalities for the mode and median of beta distributions leads to the following result.

\begin{Theorem}
\label{thm:Beta.expo}
Let $p_\ell := (a+1)/(a+b)$ and $p_r := (a-1)/(a+b)$. Then,
\[
	B_{a,b}(x) \ \le \ \begin{cases}
		x^a
			& \text{for} \ \ x \in [0,1] \ \ \text{if} \ \ b \le 1 , \\[0.5ex]
		\displaystyle
		\Bigl( \frac{x}{p_\ell} \Bigr)^a \Bigl( \frac{1 - x}{1 - p_\ell} \Bigr)^{a(b-1)/(a+1)}
			& \text{for} \ \ x \in [0,p_\ell] \ \ \text{if} \ \ b > 1 , \\[1.5ex]
		\displaystyle
		2_{}^{-1} \Bigl( \frac{x}{p} \Bigr)^a \Bigl( \frac{1 - x}{1 - p} \Bigr)^{a(b-1)/(a+1)}
			& \text{for} \ \ x \in [0,p] \ \ \text{if} \ \ a \ge b > 1 ,
	\end{cases}
\]
and
\[
	\bar{B}_{a,b}(x) \ \le \ \begin{cases}
		(1-x)^b
			& \text{for} \ \ x \in [0,1] \ \ \text{if} \ \ a \le 1 , \\[0.5ex]
		\displaystyle
		\Bigl( \frac{x}{p_r} \Bigr)^{b(a-1)/(b+1)} \Bigl( \frac{1 - x}{1 - p_r} \Bigr)^b
			& \text{for} \ \ x \in [p_r,1] \ \ \text{if} \ \ a > 1 , \\[1.5ex]
		\displaystyle
		2_{}^{-1} \Bigl( \frac{x}{p} \Bigr)^{b(a-1)/(b+1)} \Bigl( \frac{1 - x}{1 - p} \Bigr)^b
			& \text{for} \ \ x \in [p,1] \ \ \text{if} \ \ b \ge a > 1 ,
	\end{cases}
\]
\end{Theorem}

\begin{figure}
\centering
\includegraphics[width=0.85\textwidth]{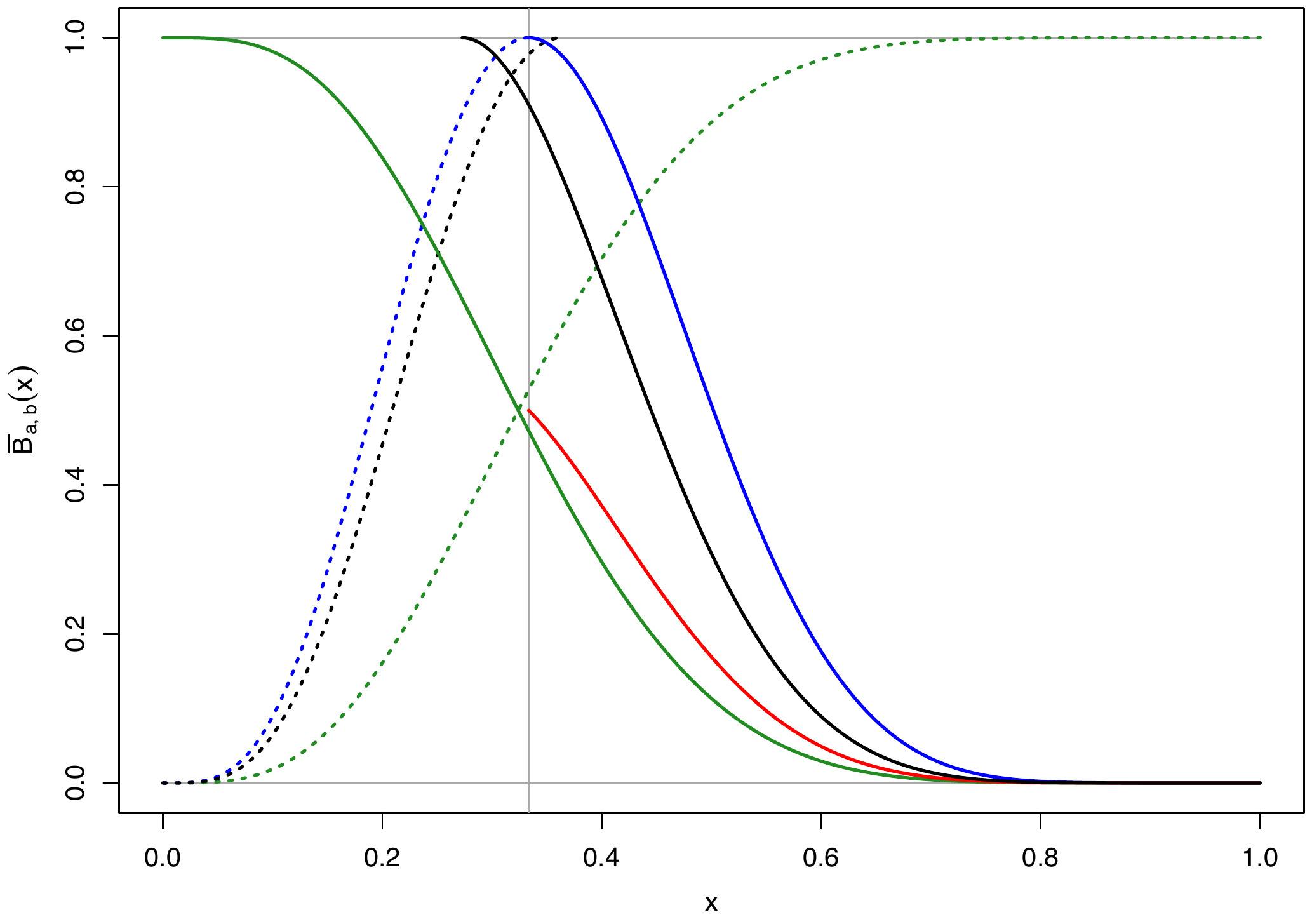}

\caption{Exponential tail inequalites for $\mathrm{Beta}(a,b)$ when $(a,b) = (4,8)$. The green line shows $\bar{B}_{a,b}$, the black and red lines are its upper bound from Theorem~\ref{thm:Beta.expo}, and the blue line is its upper bound from \eqref{ineq:Duembgen}. One also sees the distribution function $B_{a,b}$ and its bounds as dotted lines.}
\label{fig:Beta.expo}
\end{figure}

Figure~\ref{fig:Beta.expo} illustrates these bounds for $\mathrm{Beta}(4,8)$. Note also that the second upper bound for $B_{a,b}(x)$ and $\bar{B}_{a,b}(x)$ may be rewritten as follows:
\begin{align}
\label{eq:Bab1}
	B_{a,b}(x) \
	&\le \ \exp \Bigl( - \frac{a(a+b)}{a+1} K(p_\ell, x) \Bigr)
		\qquad\text{for} \ x \in [0,p_\ell] \ \ \text{if} \ \ b \ge 1 , \\
\label{eq:Bab2}
	\bar{B}_{a,b}(x) \
	&\le \ \exp \Bigl( - \frac{b(a+b)}{b+1} K(p_r, x) \Bigr)
		\qquad\text{for} \ x \in [p_r,1] \ \ \text{if} \ \ a \ge 1 .
\end{align}
This follows essentially from the fact that $c(a,b) = a(b-1)/(a+1)$ satisfies
\[
	a + c(a,b) \ = \ \frac{a(a+b)}{a+1}
	\quad\text{and}\quad
	\frac{a}{a + c(a,b)} \ = \ p_\ell ,
\]
while $c(b,a) = b(a-1)/(b+1)$ satisfies
\[
	c(b,a) + b \ = \ \frac{b(a+b)}{b+1}
	\quad\text{and}\quad
	\frac{c(b,a)}{c(b,a) + b} \ = \ p_r .
\]

\begin{Remark}[Comparison with \eqref{ineq:Duembgen} and Version~8]
For symmetry reasons, we only compare the first two upper bounds for $B_{a,b}(x)$ with the bounds in \eqref{ineq:Duembgen} and in Version~8 (Theorem~5). First of all, if $b \le 1$, then
\[
	x^a \ < \ \Bigl( \frac{x}{p} \Bigr)^a \Bigl( \frac{1 - x}{1 - p} \Bigr)^b
	\quad\text{for} \ 0 < x \le p .
\]
If $b > 1$, all bounds in question can be formulated in terms of a parameter $c \in (0,b]$ as follows:
\[
	B_{a,b}(x) \ \le \ \Bigl( \frac{x}{p(c)} \Bigr)^a \Bigl( \frac{1 - x}{1 - p(c)} \Bigr)^c
	\quad\text{for} \ 0 \le x \le p(c) ,
\]
where
\[
	p(c) \ := \ \frac{a}{a+c}
	\quad\text{and}\quad
	c \ = \ \begin{cases}
		a (b-1) /(a+1) & \text{here} , \\
		b-1 & \text{in Version~8}, \\
		b & \text{in \eqref{ineq:Duembgen}} .
	\end{cases}
\]
Note that $(b-1) a/(a+1) < b-1 < b$. Moreover, $p(c)$ is strictly decreasing in $c > 0$, and for any fixed $c_o > 0$ and $x \in (0,p(c_o)]$,
\[
	\frac{d}{dc} \log \biggl[ \Bigl( \frac{x}{p(c)} \Bigr)^a
		\Bigl( \frac{1-x}{1-p(c)} \Bigr)^c \biggr]
	= \ \log \Bigl( \frac{1 - x}{1 - p(c)} \Bigr) \
	> \ 0
\]
for $0 < c < c_o$. This implies that the present bound for $B_{a,b}(x)$ is strictly stronger than the bound of Version~8 (Theorem~5), and that the latter bound is strictly stronger than the bound in \eqref{ineq:Duembgen}.
\end{Remark}

\begin{Remark}
Note that in case of $a \ge b > 1$,
\begin{align*}
	\log \biggl[ \Bigl( \frac{x}{p_\ell} \Bigr)^a
		& \Bigl( \frac{1 - x}{1 - p_\ell} \Bigr)^{a(b-1)/(a+1)} \biggr]
	- \log \biggl[ 2_{}^{-1} \Bigl( \frac{x}{p} \Bigr)^a
		\Bigl( \frac{1 - x}{1 - p} \Bigr)^{a(b-1)/(a+1)} \biggr] \\
	&= \ \log(2) + a \log \Bigl( \frac{a}{a+1} \Bigr)
		+ \frac{a(b-1)}{a+1} \log \Bigl( \frac{b}{b-1} \Bigr) \\
	&= \ \log(2) - a \log(1 + 1/a) + \frac{a(b-1)}{a+1} \log \Bigl( \frac{b}{b-1} \Bigr) \\
	&\ge \ \log(2) - 1 + \frac{b(b-1)}{b+1} \log \Bigl( \frac{b}{b-1} \Bigr) .
\end{align*}
The right-hand side is monotone increasing in $b > 1$ and positive for $b \ge 1.5$. Thus, for $x \in [0,p]$, the third upper bound for $B_{a,b}(x)$ is smaller than the second one in case of $a \ge b \ge 1.5$.

Analogously, for $x \in [p,1]$, the third upper bound for $\bar{B}_{a,b}(x)$ is certainly better than the second one in case of $b \ge a \ge 1.5$.
\end{Remark}

The inequalities \eqref{eq:Bab1} and \eqref{eq:Bab2} imply Bernstein and Hoeffding type exponential inequalities. It follows from D\"umbgen and Wellner (2022, Lemma~S.12) and the well-known inequality $z(1-z) \le 1/4$ for $z \in \R$, that
\begin{equation}
\label{ineq:K}
	K(q,x) \
	\ge \ \frac{(x - q)^2}{2 (q/3 + 2x/3)(1 - q/3 - 2x/3)} \
	\ge \ 2 (x - q)^2
\end{equation}
for $q,x \in [0,1]$, where $K(0,x) := -\log(1 - x)$ and $K(1,x) := -\log(x)$. This leads to the following inequalities:

\begin{Corollary}
\label{cor:Beta.expo}
If $b \ge 1$, then for $x \in [0,p_\ell]$,
\begin{align*}
	B_{a,b}(x) \
	&\le \ \exp \Bigl( - \frac{a(a+b)}{a+1}
		\frac{(x - p_\ell)^2}{2(p_\ell/3 + 2x/3)(1 - p_\ell/3 - 2x/3)} \Bigr) \\
	&\le \ \exp \Bigl( - 2 \frac{a(a+b)}{a+1} (x - p_\ell)^2 \Bigr) .
\end{align*}
If $a \ge 1$, then for $x \in [p_r,1]$,
\begin{align*}
	\bar{B}_{a,b}(x) \
	&\le \ \exp \Bigl( - \frac{b(a+b)}{b+1}
		\frac{(x - p_r)^2}{2(p_r/3 + 2x/3)(1 - p_r/3 - 2x/3)} \Bigr) \\
	&\le \ \exp \Bigl( - 2 \frac{b(a+b)}{b+1} (x-p_r)^2 \Bigr) .
\end{align*}
\end{Corollary}

Further tail and concentration inequalities for the Beta distribution have been derived by Marchal and Arbel (2017) and Skorski (2023). Marchal and Arbel (2017) prove that $\mathrm{Beta}(a,b)$ is subgaussian with a variance parameter that is the solution of an equation involving hypergeometric functions. An analytic upper bound for the variance parameter is $(4(a+b+1))^{-1}$, which implies the tail inequalities
\begin{align*}
	\exp\bigl(-2(a+b+1)(x-p)^2\bigr) \ \ge \ \begin{cases}
		B_{a,b}(x) & \text{if} \ x \le p , \\
		\bar{B}_{a,b}(x) & \text{if} \ x \ge p.
	\end{cases}
\end{align*}
These bounds are weaker than the one-sided bounds in Corollary~\ref{cor:Beta.expo}. For the left tails, the difference
\[
	\frac{a(a+b)}{a+1}(x - p_{\ell})^2 - (a+b+1)(x - p)^2
\]
is strictly concave in $x$ with value $ab/(a+b)^2$ for $x = 0$ and $a/[(a+1)(a+b)]$ for $x = p$. Analogously, for the right tails, the difference
\[
	\frac{b(a+b)}{b+1} (x-p_r)^2 - (a+b+1)(x-p)^2
\]
is strictly concave in $x$ with value $b/[(b+1)(a+b)]$ for $x = p$ and $ab/(a+b)^2$ for $x = 1$.

Skorski (2023) derives a Bernstein type inequality. With the parameters
\[
	\sigma_{a,b} \ := \ \sqrt{\frac{p(1 - p)}{a+b+1}} ,
	\quad
	\gamma_{a,b}^{} \ := \ \frac{(2/3)(1-2p)^+}{a+b+2} ,
\]
he shows that for $X \sim \mathrm{Beta}(a,b)$ and $\veps \ge 0$,
\begin{equation}
\label{ineq:Skorski}
	P(X - p \ge \veps) \
	\le \ \exp \Bigl( - \frac{\veps^2}{2(\sigma_{a,b}^2 + \gamma_{a,b}^{} \veps)} \Bigr) .
\end{equation}
By means of Theorem~\ref{thm:Beta.expo} one can derive a similar inequality:

\begin{Corollary}
\label{cor:Bernstein1}
For arbitrary $a, b \ge 1$ and $\veps \ge 0$,
\[
	P(X - p \ge \veps) \
	\le \ \exp \Bigl( - \frac{(a+b+1 + a/b) \veps^2}{2[p(1 - p) + (2/3)(1 - 2p)^+ \veps]} \Bigr) .
\]
\end{Corollary}

For $p \ge 1/2$, i.e.\ $a \ge b$, the factor $a+b+1+a/b$ is at least $a+b+2$, and the bound in Corollary~\ref{cor:Bernstein1} is stronger than the one of Skorski~(2023). For sufficiently small $p < 1/2$ and sufficiently large $\veps$, Skorski's~(2023) bound \eqref{ineq:Skorski} can be a bit stronger than the one in Corollary~\ref{cor:Bernstein1}. But the next result shows that our upper bound from Theorem~\ref{thm:Segura1}, combined with \eqref{ineq:K}, implies even stronger inequalities as soon as $\veps$ is a sufficiently large multiple of $\sigma_{a,b} = \mathrm{Std}(X)$.

\begin{Corollary}
\label{cor:Bernstein2}
For any fixed $d > 0$, there exists a constant $c(d) > 0$ such that uniformly in $a,b \ge 1$,
\[
	P(X - p \ge \veps) \
	\le \ \exp \Bigl( - \frac{(a+b+d) \veps^2}{2(p + 2\veps/3)(1 - p - 2\veps/3)^+} \Bigr) ,
\]
whenever $\veps \ge c(d) \sigma_{a,b}$.
\end{Corollary}

The proof of Corollary~\ref{cor:Bernstein2} yields $c(d) = \sqrt{2/\pi} \exp(9d/4)$, but numerical experiments indicate that this is rather conservative. Note that for $\veps \in (0,1-p]$,
\begin{align*}
	\frac{(a+b+d) \veps^2}{2(p + 2\veps/3)(1 - p - 2\veps/3)} \
	&> \ \frac{(a+b+d) \veps^2}{2[p(1-p) + (2/3) (1 - 2p) \veps]} \\
	&\ge \ \frac{\veps^2}{2(\sigma_{a,b}^2 + \gamma_{a,b}^{} \veps)}
	\quad\text{if} \ \begin{cases}
		p \ge 1/2 \ \text{and} \ d \ge 1 , \\
		p < 1/2 \ \text{and} \ d \ge 2 .
	\end{cases}
\end{align*}

\section{Gaussian tail inequalities}
\label{sec:Gaussian.tail.inequalities}

Now suppose that $a,b > 1$. With $p_o := (a-1)/(a+b-2) \in (0,1)$, the density $\beta_{a,b}$ may be written as
\[
	\log \beta_{a,b}(x) \
	= \ \log \beta_{a,b}(p_o) - (a+b-2) K(p_o,x) ,
\]
whereas the probability density $\phi_{p_o,\sigma}$ of $\mathcal{N}(p_o,\sigma^2)$ with $\sigma := (4(a+b-2))^{-1/2}$ satisfies
\[
	\log \phi_{p_o,\sigma}(x) \
	= \ \log \phi_{p_o,\sigma}(p_o) - 2 (a+b-2) (x - p_o)^2 .
\]
In particular, $\rho := \log(\beta_{a,b}/\phi_{p_o,\sigma})$ satisfies
\[
	\rho'(x) \ = \ (a+b-2) (x - p_o) \bigl( 4 - 1/[x(1-x)] \bigr) ,
\]
and since $x(1-x) \le 1/4$, $\rho(x)$ is monotone decreasing in $x \ge p_o$ and monotone increasing in $x \le p_o$, where $\beta_{a,b} := 0$ on $\R \setminus (0,1)$. Consequently, for $x \ge p_o$,
\begin{align*}
	\bar{B}_{a,b}(x) \
	&\le \ \frac{\bar{B}_{a,b}(x)}{\bar{B}_{a,b}(p_o)} \\
	&= \ \frac{\int_x^\infty e^{\rho(t)} \phi_{p_o,\sigma}(t) \, \d t}
		{\int_{p_o}^x e^{\rho(t)} \phi_{p_o,\sigma}(t) \, \d t
			+ \int_x^\infty e^{\rho(t)} \phi_{p_o,\sigma}(t) \, \d t} \\
	&\le \ \frac{e^{\rho(x)} \int_x^\infty \phi_{p_o,\sigma}(t) \, \d t}
		{e^{\rho(x)} \int_{p_o}^x \phi_{p_o,\sigma}(t) \, \d t
			+ e^{\rho(x)} \int_x^\infty \phi_{p_o, \sigma}(t) \, \d t} \\
	&= \ \frac{\mathcal{N}(p_o, \sigma^2)([x,\infty))}
		{\mathcal{N}(p_o,\sigma^2)([p_o,\infty))} \\
	&= \ 2 \Phi \bigl( - 2 \sqrt{a+b-2} (x - p_o) \bigr) .
\end{align*}
Analogous arguments apply to $B_{a,b}(x)$ for $x \le p_o$, and we obtain the following bounds.

\begin{Lemma}
\label{lem:Beta.Phi}
For $a, b > 1$ and $p_o = (a-1)/(a+b-2)$,
\begin{align*}
	\bar{B}_{a,b}(x) \
	&\le \ 2 \Phi \bigl( - 2 \sqrt{a+b-2} (x - p_o) \bigr)
		\quad \text{for} \ \, x \ge p_o , \\
	B_{a,b}(x) \
	&\le \ 2 \Phi \bigl( 2 \sqrt{a+b-2} (x - p_o) \bigr)
		\quad \text{for} \ \, x \le p_o .
\end{align*}
\end{Lemma}

\section{Gaussian approximation of $\mathrm{Beta}(a,a)$}
\label{sec:Gaussian.approximation.Beta(a,a)}

Inspired by D{\"u}mbgen et al.~(2021), we want to compare the densities $\beta_{a,a}$ with the density $\phi_{1/2,\sigma}$ of $\mathcal{N}(1/2, \sigma^2)$ for various choices of $\sigma > 0$, where $a > 1$. Precisely, we want to determine
\[
	R(\sigma) \ := \ \max_{x \in (0,1)} \frac{\beta_{a,a}}{\phi_{1/2,\sigma}}(x) ,
\]
because for arbitrary Borel sets $S \subset \R$,
\[
	\mathrm{Beta}(a,a)(S) \
	\le \ R(\sigma) \, \mathcal{N}(1/2,\sigma^2)(S)
\]
and
\[
	\bigl| \mathrm{Beta}(a,a)(S) - \mathcal{N}(1/2,\sigma^2)(S) \bigr| \
	\le \ 1 - R(\sigma)^{-1} .
\]
Moreover, we want to find $\sigma > 0$ such that this quantity is minimal.

To determine $R(\sigma)$, note first that for fixed $a$ and $\sigma$,
\begin{align*}
	\log \frac{\beta_{a,a}}{\phi_{1/2,\sigma}}(x) \
	&= \ \log \sqrt{2\pi\sigma^2} - \log B(a,a) + \frac{(x - 1/2)^2}{2\sigma^2}
		+ (a-1) \log(x(1-x)) \\
	&= \ \log \sqrt{2\pi\sigma^2} - \log B(a,a) + \frac{(x - 1/2)^2}{2\sigma^2}
		+ (a-1) \log \bigl( 1/4 - (x - 1/2)^2 \bigr) \\
	&= \ \mathrm{const}(a,\sigma) + \frac{y}{8\sigma^2} + (a-1) \log(1 - y) ,
\end{align*}
where $y := (2x - 1)^2 \in [0,1)$. Since
\[
	\frac{d}{dy} \Bigl( \frac{y}{8\sigma^2} + (a-1) \log(1 - y) \Bigr) \
	= \ \frac{1}{8\sigma^2} - \frac{a-1}{1-y} ,
\]
the maximum of $\log(\beta_{a,a}/\phi_{1/2,\sigma})$ is attained at $x \in (0,1)$ such that $y = \bigl( 1 - 8\sigma^2(a-1) \bigr)^+$, and the resulting value of $\log R(\sigma)$ is
\begin{align*}
	\log R(\sigma) \
	= \ &\log \sqrt{2\pi} - \log B(a,a) + (a-1) \log(1/4) \\
	&+ \ \log(\sigma^2)/2
		+ \bigl( (8\sigma^2)^{-1} - a + 1 \bigr)^+
		+ (a-1) \log \min \{ 8 \sigma^2 (a-1), 1 \} \\
	= \ &\log \sqrt{2\pi} - \log B(a,a) - (2a-1/2) \log(2) \\
	&+ \ \log(8\sigma^2)/2
		+ \bigl( (8\sigma^2)^{-1} - a + 1 \bigr)^+
		+ (a-1) \log \min \{ 8 \sigma^2 (a-1), 1 \} .
\end{align*}
This is strictly monotone increasing in $8\sigma^2 \ge (a-1)^{-1}$, so we restrict our attention to values $\sigma$ in $\bigl( 0,(8(a-1))^{-1/2} \bigr]$. Then,
\begin{align}
\label{eq:R.sigma.1}
	\log R(\sigma) \
	= \ &\log \sqrt{2\pi} - \log B(a,a) - (2a-1/2) \log(2) \\
\nonumber
	&+ \ (8\sigma^2)^{-1} + (a - 1/2) \log(8\sigma^2)
		- a + 1 + (a-1) \log(a-1) .
\end{align}
Note also the Stirling type approximation
\begin{equation}
\label{eq:Stirling}
	\log \Gamma(y) \ = \ \log \sqrt{2\pi} + (y - 1/2) \log(y) - y + r(y) ,
\end{equation}
where $r(y)$ is strictly decreasing in $y > 0$ with $(12y + 1)^{-1} < r(y) < (12 y)^{-1}$ (cf.\ D{\"u}mbgen et al.\ 2021, Lemma~10). Consequently,
\begin{align*}
	\log\sqrt{2\pi} - \log B(a,a) \
	&= \ \log\sqrt{2\pi} + \log \Gamma(2a) - 2 \log \Gamma(a) \\
	&= \ (2a - 1/2) \log(2a) - 2(a - 1/2) \log(a) + r(2a) - 2 r(a) \\
	&= \ (2a - 1/2) \log(2) + \log(a)/2 + \tilde{r}(a) ,
\end{align*}
with $\tilde{r}(a) := r(2a) - 2 r(a)$. This leads to
\begin{align}
\label{eq:R.sigma.2}
	\log R(\sigma) \
	= \ &\tilde{r}(a) + \log(a)/2 \\
\nonumber
	&+ \ (8\sigma^2)^{-1} + (a - 1/2) \log(8\sigma^2)
		- a + 1 + (a-1) \log(a-1) .
\end{align}
For the particular choice of $\sigma$, there are at least three possibilities:\\[0.5ex]
\textbf{Moment matching.} \ A first candidate for $\sigma$ would be the standard deviation of $\mathrm{Beta}(a,a)$,
\[
	\sigma_1(a) \ := \ (8(a + 1/2))^{-1/2} .
\]
\textbf{Local density matching.} \ Since $\log \beta_{a,a}(x) - \log \beta_{a,a}(1/2)$ equals $- 4(a-1) (x - 1/2)^2$ plus a remainder of order $O((x - 1/2)^4)$ as $x \to 1/2$, another natural choice would be
\[
	\sigma_2(a) \ := \ (8(a-1))^{-1/2} .
\]
\textbf{Minimizing $R(\sigma)$.} \ Note that $\log R(\sigma) = \mathrm{const}(a) + (a - 1/2) \log(8\sigma^2) + (8\sigma^2)^{-1}$. Since
\[
	\frac{d}{dy} \bigl( (a - 1/2) \log(y) + y^{-1} \bigr) \
	= \ \frac{a - 1/2}{y} - \frac{1}{y^2} \
	= \ \frac{(a - 1/2)(y - (a - 1/2)^{-1})}{y^2} ,
\]
the optimal value of $\sigma$ equals
\[
	\sigma_3(a) \ := \ (8 (a - 1/2))^{-1/2} .
\]

\begin{figure}
\centering
\includegraphics[width=0.85\textwidth]{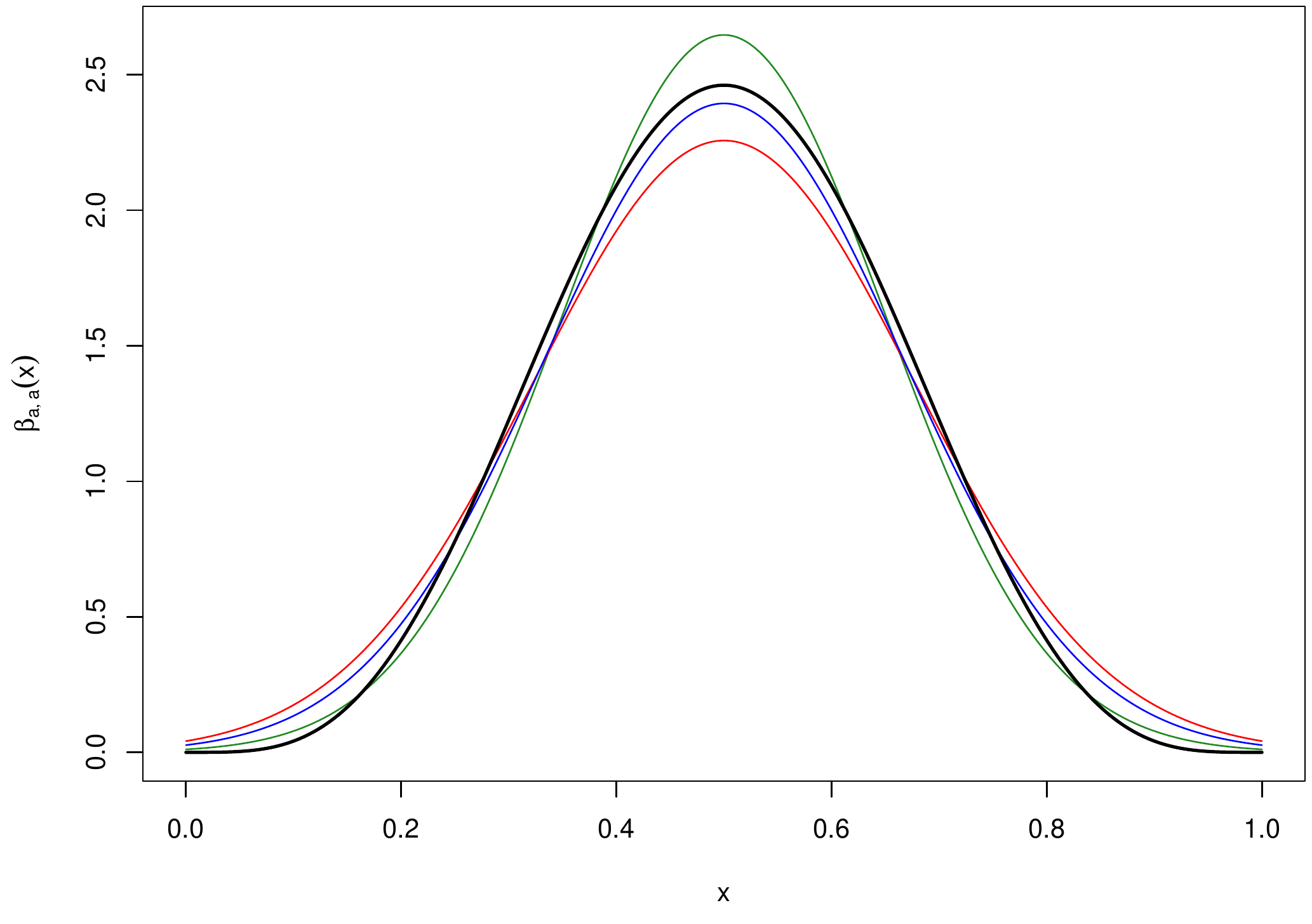}
\caption{The density $\beta_{a,a}$ (black) for $a = 5$ and its Gaussian approximation $\phi_{1/2,\sigma}$ for $\sigma = (8(a+1/2))^{-1}$ (green), $\sigma = (8(a-1))^{-1/2}$ (red) and $\sigma = (8(a - 1/2))^{-1/2}$ (blue).}
\label{fig:Beta.Gauss.A}
\end{figure}

\begin{figure}
\centering
\includegraphics[width=0.85\textwidth]{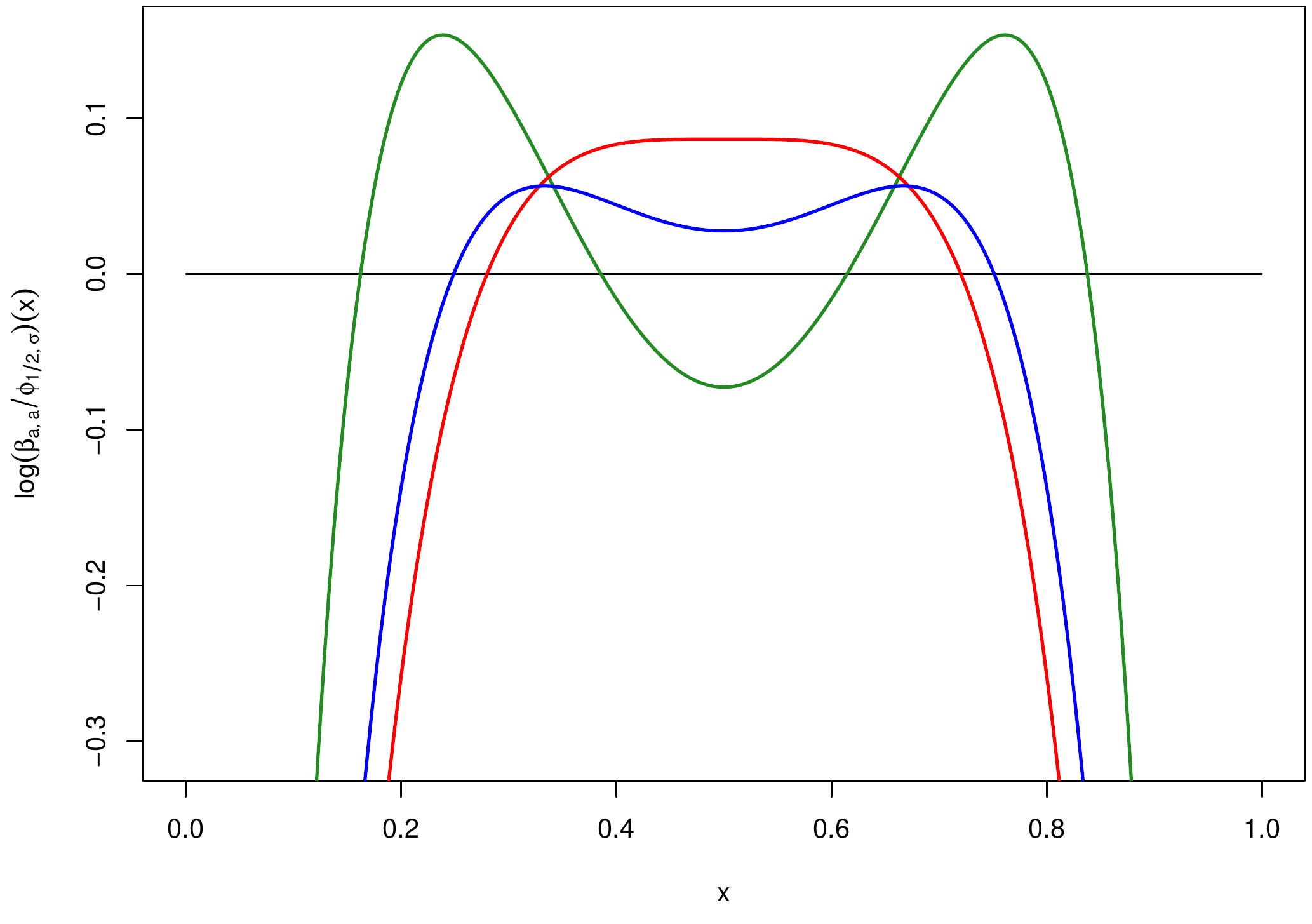}
\caption{The log-density ratios $\log(\beta_{a,a}/\phi_{1/2,\sigma})$ for $a = 5$, where $\sigma$ equals $\sigma_1(a)$ (green), $\sigma_2(a)$ (red) or $\sigma_3(a)$ (blue).}
\label{fig:Beta.Gauss.B}
\end{figure}

\paragraph{Numerical example.}
Figure~\ref{fig:Beta.Gauss.A} shows for $a = 5$ the beta density $\beta_{a,a}$ and the Gaussian approximations $\phi_{1/2,\sigma}$, where $\sigma = \sigma_1(a), \sigma_2(a), \sigma_3(a)$. Figure~\ref{fig:Beta.Gauss.B} depicts the corresponding log-density ratios $\log(\beta_{a,a}/\phi_{1/2,\sigma})$. The values of $R(\sigma)$, rounded to four digits, are $R(\sigma_1(a)) = 1.1660$, $R(\sigma_2(a)) = 1.0905$ and $R(\sigma_3(a)) = 1.0582$.

Our specific values $\sigma_j(a)$ are of the type $\sigma = (8(a + \delta))^{-1/2}$ for some $\delta \ge -1$. The next lemma provides two important properties of the resulting value $\log R(\sigma)$.

\begin{Lemma}
\label{lem:R(sigma.delta)}
Let $\sigma(a) := (8(a+\delta))^{-1/2}$ for $a > 1$ with a fixed number $\delta \ge -1$. Then $\log R(\sigma(a))$ is strictly decreasing in $a > 1$, and
\[
	\log R(\sigma) \ = \ \frac{\delta(\delta+1) + 3/4}{2a} + O(a^{-2}) .
\]
\end{Lemma} 

For our specific standard deviations $\sigma_j(a)$ we obtain the limits
\[
	\lim_{a\to \infty} \, a \log R(\sigma_j(a)) \
	= \ \begin{cases}
		3/4 & \text{if} \ j = 1 , \\
		3/8 & \text{if} \ j = 2 , \\
		1/4 & \text{if} \ j = 3 .
	\end{cases}
\]

\begin{Remark}
Similarly as in Section~\ref{sec:Gaussian.tail.inequalities}, we may conclude that for arbitrary $x > 1/2$ and $\sigma = (8(a + \delta))^{-1/2}$,
\[
	\bar{B}_{a,a}(x) \
	\le \ R(\sigma) \Phi \Bigl( - \frac{x - 1/2}{\sigma} \Bigr) \
	\le \ \frac{R(\sigma)}{2} \exp \bigl( - 4 (a + \delta) (x - 1/2)^2 \bigr) .
\]
Even the latter bound is stronger than the bound $\exp \bigl( - 4 (a + 1/2)(x - 1/2)^2 \bigr)$ by Marchal and Arbel (2017), as soon as $\delta \ge 0.5$ and $R(\sigma) \le 2$. For $\delta = 0.5$, this is the case for $a \ge 1.4$, and for $\delta = 1$, we only need $a \ge 1.9$.
\end{Remark}

\section{Proofs}
\label{sec:Proofs}

\begin{proof}[\bf Proof of Theorem~\ref{thm:Segura1}]
Let $Q$ be a continuous function on some interval $[0,x_o) \subset [0,1)$. Viewing $Q$ as a bound of $Q_{a,b}$ on $[0,x_o)$, $H(x) = x^a Q(x) /[a B(a,b)]$ is a bound for $B_{a,b}(x)$ on $[0,x_o)$ with $H(0) = 0$. If $Q$ is differentiable on $(0,x_o)$, then elementary calculus reveals that
\[
	H'(x) \ = \ \beta_{a,b}(x) J(x)
	\quad\text{with}\quad
	J(x) \ := \ \frac{Q(x) + Q'(x) x/a}{(1 - x)^{b-1}} .
\]
If we can show that $J \ge 1$ or $J \le 1$ on $(0,x_o)$, we may conclude that $Q_{a,b} \le Q$ or $Q_{a,b} \ge Q$, respectively, on $[0,x_o)$.

On the one hand, let $Q(x) := (1 - x)^b [1 + cx/(1-x)] = (1 - x)^b + c x (1 - x)^{b-1}$ for some $c > 0$ and $x \in [0,1)$. Then $Q'(x) = (c-b)(1-x)^{b-1} + (1-b)c x (1 - x)^{b-2}$, and elementary calculations lead to the formula
\begin{align*}
	J(x) \
	&= \ 1 + (1 + a^{-1}) (c - c_{a,b}) x + \frac{(1-b)c x^2}{a(1 - x)} .	
\end{align*}
If we choose $c = c_{a,b}$, then
\[
	J(x) \
	= \ 1 + \frac{(1-b)c_{a,b} x^2}{a(1 - x)} \
	\begin{cases}
		\ge \ 1 & \text{if} \ b \le 1 , \\
		\le \ 1 & \text{if} \ b \ge 1 .
	\end{cases}	
\]
This proves the bounds in terms of $Q_{a,b}^{[S,1]}$.

On the other hand, let $Q(x) := (1 - x)^b /(1 - cx)$ for some $c > 0$ and $0 \le x < x_o := \min\{c^{-1},1\}$. For $0 < x < x_o$,
\[
	Q'(x) \ = \ - \frac{b(1-x)^{b-1}}{1-cx} + \frac{c(1-x)^b}{(1 - cx)^2} ,
\]
and elementary calculations lead to
\begin{align*}
	J(x) \
	&= \ 1 + (1 + a^{-1}) (c - c_{a,b}) \frac{x}{1 - cx} + \frac{c(c-1)x^2}{a(1 - cx)^2} .
\end{align*}
If we choose $c = c_{a,b}$, then
\[
	J(x) \
	= \ 1 + \frac{c_{a,b}(c_{a,b}-1)x^2}{a(1 - c_{a,b}x)^2} \
	\begin{cases}
		\le \ 1 & \text{if} \ c_{a,b} \le 1, \ \text{i.e.} \ b \le 1 , \\
		\ge \ 1 & \text{if} \ c_{a,b} \ge 1, \ \text{i.e.} \ b \ge 1 .
	\end{cases}		
\]
Hence, $Q$ is a lower or upper bound for $Q_{a,b}$ if $b \le 1$ or $b \ge 1$, respectively. But this bound can be refined further. Note that for $x \in [0,x_o)$,
\[
	Q'(x) \ = \ \frac{(b-1)c(1 - x)^{b-1}}{(1 - cx)^2} (x - p) ,
\]
and $p < x_o$. If $b \le 1$, then $Q$ is decreasing on $[p,1)$ whereas $Q_{a,b}$ is increasing on $[0,1)$. Thus we may conclude that
\[
	Q_{a,b}(x) \ \ge \ Q_{a,b}(p) \ \ge \ Q(p) \ = \ (a+1) (1 - p)^b
	\quad\text{for} \ x \ge p .
\]
Likewise, if $b \ge 1$, then $Q$ is increasing on $[p,x_o)$, whereas $Q_{a,b}$ is decreasing on $[0,1)$. Thus we may conclude that
\[
	Q_{a,b}(x) \ \le \ Q_{a,b}(p) \le \ Q(p) \ = \ (a+1) (1 - p)^b
	\quad\text{for} \ x \ge p .
\]
This proves the bounds in terms of $Q_{a,b}^{[S,2]}$.
\end{proof}

\begin{proof}[\bf Proof of Theorem~\ref{thm:Segura2}]
Note first that for any $x \in [0,1)$, the integral $a \int_0^1 y^{a-1} (1 - xy)^{b-1} \, \d y = Q_{a,b}(x)$ is well-defined for any $a > 0$ and $b \in \R$. We may also write
\[
	Q_{a,b}(x) \
	= \ \Ex[(1 - x Y)^{b-1}]
\]
with $Y \sim \mathrm{Beta}(a,1)$. Since $\Ex(Y) = a/(a+1)$ and
\[
	\frac{\d^2}{\d y^2} (1 - xy)^{b-1} \ = \ (b-1)(b-2) x^2 (1 - xy)^{b-3} \
	\begin{cases}
		\ge 0 & \text{if} \ b \in (-\infty,1] \cup [2,\infty) , \\
		\le 0 & \text{if} \ b \in [1,2] ,
	\end{cases}
\]
it follows from Jensen's inequality that
\[
	Q_{a,b}(x) \ \begin{cases}
		\ge \ Q_{a,b}^{[1]}(x) & \text{if} \ b \in (-\infty,1] \cup [2,\infty) , \\
		\le \ Q_{a,b}^{[1]}(x) & \text{if} \ b \in [1,2] ,
	\end{cases}
\]
where
\[
	Q_{a,b}^{[1]}(x) \ := \ \bigl( 1 - x \Ex(Y) \bigr)^{b-1} \
	= \ \Bigl( 1 - \frac{ax}{a+1} \Bigr)^{b-1} .
\]

Secondly, it follows from partial integration that
\[
	Q_{a,b}(x) \ = \ (1 - x)^{b-1} + \frac{b-1}{a+1} x Q_{a+1,b-1}(x) .
\]
Consequently
\[
	Q_{a,b}(x) \ \begin{cases}
		\ge \ Q_{a,b}^{[2]}(x) & \text{if} \ b \in [1,2] \cup [3,\infty) , \\
		\le \ Q_{a,b}^{[2]}(x) & \text{if} \ b \in (-\infty,1] \cup [2,3] ,
	\end{cases}
\]
where
\begin{align*}
	Q_{a,b}^{[2]}(x) \
	:= \ &(1 - x)^{b-1} + \frac{b-1}{a+1} x Q_{a+1,b-1}^{[1]}(x) \\
	= \ &(1 - x)^{b-1} + \frac{(b-1)x}{a+1} \Bigl( 1 - \frac{(a+1)x}{a+2} \Bigr)^{b-2} .
\end{align*}

Thirdly, note that
\[
	\frac{\d}{\d x} Q_{a,b}(x) \
	= \ - a (b-1) \int_0^1 y^a (1 - xy)^{b-2} \, \d y \
	= \ - \frac{a(b-1)}{a+1} Q_{a+1,b-1}(x) .
\]
Since $Q_{a,b}(0) = 0$, this implies that
\[
	Q_{a,b}(x) \ = \ 1 - \frac{a(b-1)}{a+1} \int_0^x Q_{a+1,b-1}(t) \, \d t .
\]
In particular,
\[
	Q_{a,b}(x) \ \begin{cases}
		\le \ Q_{a,b}^{[3]}(x) & \text{if} \ b \in [1,2] \cup [3,\infty) , \\
		\ge \ Q_{a,b}^{[3]}(x) & \text{if} \ b \in (-\infty,1] \cup [2,3] ,
	\end{cases}
\]
where
\begin{align*}
	Q_{a,b}^{[3]}(x) \
	:= \ &1 - \frac{a(b-1)}{a+1} \int_0^x Q_{a+1,b-1}^{[1]}(t) \ \d t \\
	= \ &1 - \frac{a(b-1)}{a+1} \int_0^x \Bigl( 1 - \frac{(a+1)t}{a+2} \Bigr)^{b-2} \, \d t \\
	= \ &1 + \frac{a(a+2)}{(a+1)^2} \Bigl( 1 - \frac{(a+1)t}{a+2} \Bigr)^{b-1} \Big|_{t=0}^x \\
	= \ &\frac{1}{(a+1)^2} + \frac{a(a+2)}{(a+1)^2} \Bigl( 1 - \frac{(a+1)x}{a+2} \Bigr)^{b-1} .
\end{align*}

Finally, to derive the bounds in terms of $Q_{a,b}^{[4]}(x)$ and $Q_{a,b}^{[5]}(x)$, we use a particular representation of a twice continuously differentiable function $f$ on $[0,1]$. Namely, for any $y \in (0,1)$,
\begin{align*}
	f(y) \
	&= \ (1 - y) f(0) + y f(1)
		- y(1-y) \int_0^1 \min \Bigl( \frac{t}{y}, \frac{1-t}{1-y} \Bigr) f''(t) \, \d t \\
	&= \ (1 - y) f(0) + y f(1) - \frac{y(1-y)}{2} \Ex f''(\Delta_y) ,	
\end{align*}
where $\Delta_y$ is a random variable with values in $(0,1)$ and density function
\[
	g_y(t) \ := \ 2 \min \Bigl( \frac{t}{y}, \frac{1-t}{1-y} \Bigr) , \quad t \in (0,1) .
\]
Now we apply this to the function $f(y) = f_x(y) := (1 - xy)^{b-1}$. Here, $f_x(0) = 1$ and $f_x(1) = (1-x)^{b-1}$, while $f_x''(u) = (b-1)(b-2) x^2 (1 - xy)^{b-3}$. Thus,
\begin{align*}
	Q_{a,b}(x) \
	= \ &a \int_0^1 \bigl( y^{a-1}(1 - y) + y^a (1 - x)^{b-1} \bigr) \, \d y \\
		&- \ \frac{a (b-1)(b-2)x^2}{2} \int_0^1 y^a(1 - y) \Ex[(1 - x \Delta_y)^{b-3}] \, \d y \\
	&= \ \frac{a(1 - x)^{b-1} + 1}{a+1}
		- \frac{a(b-1)(b-2)x^2}{2} \int_0^1 y^a(1 - y) \Ex[(1 - x \Delta_y)^{b-3}] \, \d y \\
	&= \ \frac{a(1 - x)^{b-1} + 1}{a+1}
		- \frac{a(b-1)(b-2)x^2}{2(a+1)(a+2)} \Ex \bigl[ \Ex[(1 - x \Delta_Y)^{b-3} \,|\, Y] \bigr] ,
\end{align*}
where $Y \sim \mathrm{Beta}(a+1,2)$ and, conditional on $Y$, $\Delta_Y$ follows the density $g_Y$. On the one hand, $\Ex(Y) = (a+1)/(a+3)$, and elementary calculations reveal that $\Ex(\Delta_y) = (1 + y)/3$. On the other hand, $z^{b-3}$ is convex or concave in $z > 0$ if $b \not\in (3,4)$ or $b \in [3,4]$, respectively. Consequently, it follows from a two-fold application of Jensen's inequality that
\[
	\Bigl( 1 - \frac{1 + \Ex(Y)}{3} x \Bigr)^{b-3} \
	= \ \Bigl( 1 - \frac{2(a+2)x}{3(a+3)} \Bigr)^{b-3} \
\]
is a lower or upper bound for $\Ex \bigl[ \Ex[(1 - x \Delta_Y)^{b-3} \,|\, Y] \bigr]$ if $b \not\in (3,4)$ or $b \in [3,4]$, respectively. Taking into account the sign of $(b-1)(b-2)$, we see that
\[
	Q_{a,b}^{[4]}(x) \ = \ \frac{a(1 - x)^{b-1} + 1}{a+1}
		- \frac{a(b-1)(b-2)x^2}{2(a+1)(a+2)} \Bigl( 1 - \frac{2(a+2)x}{3(a+3)} \Bigr)^{b-3}
\]
satisfies the following inequalities:
\[
	Q_{a,b}(x) \ \begin{cases}
		\le \ Q_{a,b}^{[4]}(x) & \text{if} \ b \in (0,1] \cup [2,3] \cup [4,\infty) , \\[0,5ex]
		\ge \ Q_{a,b}^{[4]}(x) & \text{if} \ b \in [1,2] \cup [3,4] .
	\end{cases}
\]
Instead of using Jensen's inequality, one can use the fact that the affine interpolant
\[
	1 - \Delta_y + \Delta_y (1 - x)^{b-3}
\]
is an upper or lower bound for $(1 - x \Delta_y)^{b-3}$ if $b \not\in (3,4)$ or $b \in [3,4]$, respectively, and
\[
	\Ex(1 - \Delta_Y + \Delta_Y (1 - x)^{b-3}) \
	= \ \frac{a + 5 + 2(a+2) (1 - x)^{b-3}}{3(a+3)} .
\]
Taking into account the sign of $(b-1)(b-2)$, we see that
\[
	Q_{a,b}^{[5]}(x) \ = \ \frac{a(1 - x)^{b-1} + 1}{a+1}
		- \frac{a(b-1)(b-2)x^2 \bigl( a + 5 + 2(a+2)(1-x)^{b-3} \bigr)}{6(a+1)(a+2)(a+3)}
\]
satisfies the following inequalities:
\[
	Q_{a,b}(x) \ \begin{cases}
		\ge \ Q_{a,b}^{[5]}(x) & \text{if} \ b \in (0,1] \cup [2,3] \cup [4,\infty) , \\[0,5ex]
		\le \ Q_{a,b}^{[5]}(x) & \text{if} \ b \in [1,2] \cup [3,4] .
	\end{cases}
\]\\[-5ex]
\end{proof}

\begin{proof}[\bf Proof of Remark~\ref{rem:Segura12}]
We start with the lower bounds for $Q_{a,b}$. For $x \in (0,1]$ and $b \in (0,1)$, the asserted inequality $Q_{a,b}^{[S,2]}(x) < Q_{a,b}^{[1]}(x)$ may be rewritten as
\[
	b \log(1 - x) - \log[1 - x + (1 - b) x/(a+1)] \ < \ (b-1) \log[1 - x + x/(a+1)] ,
\]
and this is equivalent to
\[
	\log[1 + (1-b)y] \ > \ (1-b) \log(1 + y)
\]
with $y := x/[(a+1)(1-x)] > 0$. But the latter inequality is a direct consequence of strict concavity of $\log(\cdot)$.

For $x \in (0,1]$ and $b \in (1,2]$, we write
\[
	Q_{a,b}^{[S,1]}(x) \
	= \ (1 - x)^{b-1} (1 - x + c_{a,b}x) \
	= \ (1 - x)^{b-1} + \frac{(b-1)x}{a+1} \, (1 - x)^{b-1} ,
\]
so the assertion that $Q_{a,b}^{[S,1]}(x) < Q_{a,b}^{[2]}(x)$ is equivalent to
\[
	(1 - x)^{b-1} \ < \ [1 - (a+1)x/(a+2)]^{b-2} .
\]
But the right-hand side is strictly larger than $[1 - (a+1) x/(a+2)]^{b-1} > (1 - x)^{b-1}$.

For $x \in (0,1]$ and $b \in [2,\infty)$, the asserted inequality $Q_{a,b}^{[S,1]}(x) < Q_{a,b}^{[1]}(x)$ is equivalent to
\[
	b \log(1 - x) + \log \Bigl( 1 + \frac{(a+b)x}{(a+1)(1-x)} \Bigr) \
	< \ (b-1) \log \Bigl( 1 - \frac{ax}{a+1} \Bigr) .
\]
This inequality is even valid for all $b > 1$. Indeed, for $b = 1$, both sides coincide, and the derivatives of the left-hand and right-hand side with respect to $b$ are equal to
\[
	\log(1 - x) + \frac{x}{a+1 + (b-1) x}
	\quad\text{and}\quad
	\log \Bigl( 1 - \frac{ax}{a+1} \Bigr) ,
\]
respectively. The former is strictly decreasing while the latter is constant in $b \ge 1$. Thus it suffices to show that
\[
	\log(1 - x) + \frac{x}{a+1} \
	\le \ \log \Bigl( 1 - \frac{ax}{a+1} \Bigr) .
\]
Writing $1 - ax/(a+1) = 1 - x + x/(a+1)$, an equivalent claim is that
\[
	\frac{x}{a+1} \ \le \ \log \Bigl( 1 + \frac{x}{(a+1)(1-x)} \Bigr) .
\]
But $\log(1+y) \ge y/(1+y)$ for $y \ge 1$, whence
\[
	\log \Bigl( 1 + \frac{x}{(a+1)(1-x)} \Bigr) \
	\ge \ \frac{x}{a+1 - ax} \
	> \ \frac{x}{a+1} .
\]

Now we consider the upper bounds for $Q_{a,b}$. For $x \in (0,1]$ and $b \in (0,1)$, the inequality $Q_{a,b}^{[S,1]}(x) > Q_{a,b}^{[2]}$ can be rewritten as
\[
	1 + \frac{(b-1)x}{a+1} \
	> \ 1 + \frac{(b-1)x}{a+1} \Bigl( 1 + \frac{x}{a+2} \Bigr)^{b-1}
		\Bigl( 1 - \frac{(a+1)x}{a+2} \Bigr)^{-1} ,
\]
and since $b-1 < 0$, this is equivalent to
\[
	(1 + y)^{1-b} (1 - (a+1)y) \ < \ 1
\]
with $y := x/(a+2) \in (0,(a+2)^{-1}]$. But the left-hand side is smaller than $(1 + y)(1-y) = 1 - y^2 < 1$.

For $x \in (0,p]$ and $b \in (1,2]$, the inequality $Q_{a,b}^{[S,2]}(x) > Q_{a,b}^{[1]}(x)$ can be rewritten as
\[
	b \log(1 - x) - \log[1 - x - (b - 1) x/(a+1)] \ > \ (b-1) \log[1 - x + x/(a+1)] ,
\]
and this is equivalent to
\[
	\log[1 - (b-1) y] + (b-1) \log(1 + y) \ < \ 0
\]
with $y := x/[(a+1)(1-x)] > 0$. But the left-hand side equals
\begin{align*}
	b &\bigl( b^{-1} \log[1 - (b-1)y] + (1 - b^{-1}) \log(1 + y) \bigr) \\
	&< \ b \log \bigl( b^{-1}[1 - (b-1)y] + (1 - b^{-1})(1 + y) \bigr) \
		= \ b \log(1) \ = \ 0
\end{align*}
by strict concavity of $\log(\cdot)$. For $x \in (p,1)$, the inequality $Q_{a,b}^{[S,2]}(x) > Q_{a,b}^{[1]}(x)$ follows from $Q_{a,b}^{[S,2]}$ being constant and $Q_{a,b}^{[1]}$ being decreasing on $[p,1)$.

Finally, for $x \in (0,p]$ and $b \in [2,3]$, the inequality $Q_{a,b}^{[S,2]}(x) > Q_{a,b}^{[2]}(x)$ can be rewritten as
\[
	(1 - y)^{-1} \ > 1 + y [1 + (a+1)y/(a+2)]^{b-2}
\]
with $y := x/[(a+1)(1-x)] > 0$. But since $(a+1)/(a+2) < 1$ and $0 \le b-2 \le 1$,
\[
	1 + y [1 + (a+1)y/(a+2)]^{b-2} \
	\le \ 1 + y(1+y) \
	= \ 1 + y + y^2 \
	< \ (1 - y)^{-1} .
\]\\[-5ex]
\end{proof}

\begin{proof}[\bf Proof of Corollary~\ref{cor:Segura}]
The Stirling approximation \eqref{eq:Stirling} implies the following asymptotic expansions:
\[
	\frac{1}{B(a,b)} \
	= \ \frac{\Gamma(a+b)}{\Gamma(a) \Gamma(b)} \
	= \ \frac{b^a (1 + a/b)^{a+b-1/2} e^{-a} (1 + o(1))}{\Gamma(a)} \
	= \ \frac{b^a (1 + o(1))}{\Gamma(a)} .
\]
Throughout this proof, asymptotic statements refer to $b \to \infty$. Consequently,
\[
	G_a(x) \
	= \ \lim_{b\to \infty} B_{a,b}(x/b) \
	= \ \lim_{b\to\infty} \frac{(x/b)^a}{a B(a,b)} Q_{a,b}(x/b) \
	= \ \frac{x^a}{a\Gamma(a)} \lim_{b\to \infty} Q_{a,b}(x/b) .
\]
Now we bound $Q_{a,b}(x/b)$ in terms of $Q_{a,b}^{[S,2]}(x/b)$ as in Theorem~\ref{thm:Segura1} or $Q_{a,b}^{[\ell]}(x/b)$, $1 \le \ell \le 5$, as in Theorem~\ref{thm:Segura2}. Elementary calculations reveal that
\begin{align*}
	Q_{a,b}^{[S,2]}(x/b) \
	&\to \ Q_a^{[S]}(x) \
		:= \ \frac{a+1}{(a+1-x)^+} \, e_{}^{-x} , \\
	Q_{a,b}^{[1]}(x/b) \
	&\to \ Q_a^{[1]}(x) \
		:= \ e_{}^{-ax/(a+1)} , \\
	Q_{a,b}^{[2]}(x/b) \
	&\to \ Q_a^{[2]}(x) \
		:= \ e_{}^{-x} + \frac{x}{a+1} \, e_{}^{-(a+1)x/(a+2)} , \\
	Q_{a,b}^{[3]}(x/b) \
	&\to \ Q_a^{[3]}(x) \
		:= \ \frac{1}{(a+1)^2} + \frac{a(a+2)}{(a+1)^2} \, e_{}^{-(a+1)x/(a+2)} , \\
	Q_{a,b}^{[4]}(x/b) \
	&\to \ Q_a^{[4]}(x) \
		:= \ \frac{ae^{-x} + 1}{a+1} - \frac{ax^2}{2(a+1)(a+2)} \, e_{}^{-2(a+2) x/[3(a+3)]} , \\
	Q_{a,b}^{[5]}(x/b) \
	&\to \ Q_a^{[5]}(x) \
		:= \ \frac{ae^{-x} + 1}{a+1} - \frac{ax^2 [a+5 + 2(a+2) e^{-x}]}{6(a+1)(a+2)(a+3)} .
\end{align*}
Hence the asserted bounds follow from the fact that for $b > 4$, $Q_{a,b}^{[1]}, Q_{a,b}^{[1]}, Q_{a,b}^{[2]}$ are lower and $Q_{a,b}^{[S,2]}, Q_{a,b}^{[3]}, Q_{a,b}^{[4]}$ are upper bounds for $Q_{a,b}$.

As to $\bar{G}_a(x)$, we write $\bar{G}_a(x) = \lim_{b \to \infty} \bar{B}_{a,b}(x/b) = \lim_{b \to \infty} B_{b,a}(1 - x/b)$ and
\[
	B_{b,a}(1 - x/b) \
	= \ \frac{(1 - x/b)^b}{b B(a,b)} Q_{b,a}(1 - x/b) \
	= \ \frac{e^{-x} b^{a-1} (1 + o(1))}{\Gamma(a)} Q_{b,a}(1 - x/b) ,
\]
so $\bar{G}_a(x)$ is $e^{-x}/\Gamma(a)$ times $\lim_{b \to \infty} b^{a-1} Q_{b,a}(1 - x/b)$. Now we bound $Q_{b,a}(1 - x/b)$ in terms of $Q_{b,a}^{[S,2]}(1 - x/b)$ as in Theorem~\ref{thm:Segura1} or $Q_{b,a}^{[1]}(1 - x/b), Q_{b,a}^{[2]}(1 - x/b)$ as in Theorem~\ref{thm:Segura2}. (The bounds in terms of $Q_{b,a}^{[\ell]}(1 - x/b)$, $3 \le \ell \le 5$, turned out to be useless in this particular context.) One can show that
\begin{align*}
	b^{a-1} Q_{b,a}^{[S,2]}(1 - x/b) \
	&\to \ \bar{Q}_a^{[S]}(x) \
		:= \ \frac{x^a}{(x - a + 1)^+} , \\
	b^{a-1} Q_{b,a}^{[1]}(1 - x/b) \
	&\to \ \bar{Q}_a^{[1]}(x) \
		:= (x + 1)^{a-1} , \\
	b^{a-1} Q_{b,a}^{[2]}(1 - x/b) \
	&\to \ \bar{Q}_a^{[2]}(x) \
		:= x^{a-1} + (a-1) (x+1)^{a-2} .
\end{align*}
Now the asserted bounds for $\bar{G}_a$ are a consequence of Theorems~\ref{thm:Segura1} and \ref{thm:Segura2} with $(b,a)$ in place of $(a,b)$.
\end{proof}

\begin{proof}[\bf Proof of Proposition~\ref{prop:Beta}]
For symmetry reasons, it suffices to show that for $x \in [0,q]$, $B_{a,b}(x)$ is not larger than
\[
	\tilde{B}(x) \ := \ B_{a,b}(q) \Bigl( \frac{x}{q} \Bigr)^a \Bigl( \frac{1 - x}{1 - q} \Bigr)^{c}
\]
with $c = a (b-1)^+/(a+1)$. Note that $\tilde{B}(0) = B(0) = 0$ and $\tilde{B}(q) = B(q) > 0$. Moreover, for some positive constants $d, \tilde{d}$ depending on $a$, $b$ and $q$,
\begin{align*}
	\tilde{B}'(x) \
	&= \ d \bigl( a x^{a-1} (1 - x)^c - c x^a (1 - x)^{c-1} \bigr) \\
	&= \ d x^{a-1} (1 - x)^{b-1} (1 - x)^{c-b} (a - (a+c) x) \\
	&= \ B_{a,b}'(x) J(x)
\end{align*}
with
\[
	J(x) \ := \ \tilde{d} (1 - x)^{c-b} (a - (a+c) x) .
\]
If we can show that $J$ is monotone decreasing on $[0,q]$, then a standard argument for measures with monotone density ratios applies: If $J(x) \ge 1$, then $\tilde{B}(x) - B_{a,b}(x) = \int_0^x (J(t) - 1) B_{a,b}'(t) \,\d t \ge 0$, whereas if $J(x) \le 1$, then $\tilde{B}(x) - B_{a,b}(x) = \int_x^q (1 - J(t)) B_{a,b}'(t) \,\d t \ge 0$.

As to the monotonicity of $J$,
\begin{align*}
	J'(x) \
	&= \ \tilde{d} \bigl[ - (c - b)(1 - x)^{c-b-1} (a - (a+c)x) - (a+c) (1 - x)^{c-b} \bigr] \\
	&= \ \tilde{d} (1 - x)^{c-b-1} \bigl[ (b-1)a - (a+1)c - (b-1-c)(a+c) x \bigr] \\
	&= \ \tilde{d} (1 - x)^{c-b-1} \cdot \begin{cases}
			(b-1)a (1 - x) & \text{if} \ b \le 1 , \\[1ex]
			\displaystyle
			- \frac{(b-1) a (a+b)}{(a+1)^2} \, x & \text{if} \ b \ge 1 ,
		\end{cases} 
\end{align*}
and this is nonpositive indeed.
\end{proof}

\begin{proof}[\bf Proof of Theorem~\ref{thm:Beta.expo}]
Again, for symmetry reasons, it suffices to derive the upper bounds for $B_{a,b}$. If $b \le 1$, then Proposition~\ref{prop:Beta} shows that $B_{a,b}(x) \le (x/q)^a$ for $x \in [0,q]$ and $q \in (0,1)$. Letting $q \to 1$ reveals that $B_{a,b}(x) \le x^a$ for $x \in [0,1]$. If $b > 1$, the maximizer of $q \mapsto q^a (1 - q)^{c(a,b)}$ equals $a/(a + c(a,b)) = (a+1)/(a+b) = p_\ell$. Thus Proposition~\ref{prop:Beta} yields the second bound for $B_{a,b}(x)$, $x \in [0,p_\ell]$. Moreover, it is well-known that in case of $a \ge b > 1$, the median of $\mathrm{Beta}(a,b)$ is at least its mean $a/(a+b)$, see Groeneveld and Meeden~(1977) or Dharmadhikari and Joag-Dev~(1983). Thus $B_{a,b}(p) \le 2^{-1}$, whence Proposition~\ref{prop:Beta} leads to the third asserted bound for $B_{a,b}(y)$.
\end{proof}

\begin{proof}[\bf Proof of Corollary~\ref{cor:Bernstein1}]
It suffices to prove the assertion for $\veps \in (0,1-p]$, because $X \le 1$ almost surely. Corollary~\ref{cor:Beta.expo} shows that for $a,b \ge 1$,
\begin{align*}
	- \log \Pr(X \ge p + \veps) \
	&= \ - \log \bar{B}_{a,b}(p + \veps) \\
	&\ge \ \frac{m (1 + 1/b)^{-1} (\veps + \delta)^2}{2 [p + (2\veps - \delta)][1 - p - (2\veps - \delta)]} ,
\end{align*}
where $m := a+b$ and $\delta := 1/m$. Thus it suffices to show that the right-hand side of the previous display is not smaller than
\[
	\frac{(a + b + 1 + a/b) \veps^2}{2 [p(1-p) + (1 - 2p)^+ 2\veps/3]} \
	= \ \frac{m(1 + 1/b) \veps^2}{2 [p(1-p) + (1 - 2p)^+ 2\veps/3]} . 
\]
But this claim is equivalent to
\begin{equation}
\label{ineq:Bernstein1A}
	\Bigl( \frac{1 + \delta/\veps}{1 + 1/b} \Bigr)^2 \
	\ge \ \frac{p(1-p) + (1 - 2p) (2\veps - \delta)/3 - (2\veps - \delta)^2/9}
		{p(1 - p) + (2/3) (1 - 2p)^+ \veps} .
\end{equation}
For $p \le 1/2$, this follows from the facts that the left-hand side of \eqref{ineq:Bernstein1A} is at least $1$, because $\delta/\veps \ge \delta/(1-p) = 1/b$, and that the right-hand side of \eqref{ineq:Bernstein1A} equals
\[
	\frac{p(1-p) + (2/3)(1 - 2p) \veps - (1 - 2p)\delta/3 - (2\veps - \delta)^2/9}
		{p(1 - p) + (2/3) (1 - 2p) \veps} \
	\le \ 1 .
\]
Now suppose that $p > 1/2$. Then the right-hand side of \eqref{ineq:Bernstein1A} equals
\[
	\frac{p(1-p) - (2p - 1) (2\veps - \delta) - (2\veps - \delta)^2/9}{p(1 - p)} \
	\ \le \ 1 - \frac{2p-1}{p(1-p)}(2\veps - \delta) ,
\]
so it suffices to verify that
\[
	\Bigl( \frac{1 + \delta/\veps}{1 + 1/b} \Bigr)^2 \
	\ge \ 1 - \frac{2p-1}{p(1-p)}(2\veps - \delta) .
\]
Since the left-hand side is at least $1$, this claim is obvious for $\veps \ge \delta/2$. To verify the inequality for $\veps \in (0,\delta/2]$, it suffices to show that
\[
	\frac{\d}{\d \veps} \biggl( \Bigl( \frac{1 + \delta/\veps}{1 + 1/b} \Bigr)^2 
		- 1 + \frac{2p-1}{p(1-p)}(2\veps - \delta) \biggr) \
	\le \ 0
\]
for $\veps \le \delta/2$. Indeed, the left-hand side equals
\begin{align*}
	- \frac{2 \delta (1 + \delta/\veps)}{\veps^2 (1 + 1/b)^2} + \frac{2(2p-1)}{p(1-p)} \
	&\le \ - \frac{24}{\delta (1 + 1/b)^2} + \frac{2(2p-1)}{p(1-p)} \\
	&= \ - \frac{24 m}{(1+1/b)^2} + 2m(1/b-1/a) \\
	&< \ - 4 m ,
\end{align*}
because $1/b \le 1$ and $1/a > 0$.
\end{proof}

\begin{proof}[\bf Proof of Corollary~\ref{cor:Bernstein2}]
If $a,b \ge 1$, it follows from Theorem~\ref{thm:Segura1} that for $x \in [p,1]$,
\begin{align*}
	- \log & \Pr(X \ge x) \\
	&= \ - \log \bar{B}_{a,b}(x) \ = \ - \log B_{b,a}(1 - x) \\
	&\ge \ \log[b B(a,b)] - b \log(1 - x) - a \log(x) + \log [1 - c_{b,a}(1 - x)] \\
	&= \ \log B(a,b) - m p \log(x) - m (1-p) \log(1 - x) + \log(1 - a + mx) + \log[b/(b+1)] \\
	&\ge \ \log B(a,b) - m p \log(x) - m (1-p) \log(1 - x) + \log(1 - a + mx) - \log(2)	
\end{align*}
with $m = a+b$. But it follows from the Stirling approximation \eqref{eq:Stirling} that
\begin{align*}
	\log B(a,b) \
	&= \ \log \Gamma(a) + \log \Gamma(b) - \log \Gamma(a+b) \\
	&\ge \ (a - 1/2) \log(a) + (b - 1/2) \log(b) - (m - 1/2) \log(m) + \log \sqrt{2\pi} \\
	&= \ mp \log(p) + m(1-p) \log(1 - p) - \log \sqrt{mp(1-p)} + \log \sqrt{2\pi} ,
\end{align*}
whence
\begin{align*}
	- \log \Pr(X \ge x) \
	&\ge \ m K(p,x) + \log[1 + m(x-p)] - \log \sqrt{mp(1-p)} + \log \sqrt{\pi/2} \\
	&\ge \ m K(p,x) + \log[(m+1)(x-p)] - \log \sqrt{mp(1-p)} + \log \sqrt{\pi/2} \\
	&\ge \ m K(p,x) + \log[(x - p)/\sigma_{a,b}] + \log\sqrt{\pi/2} .
\end{align*}
Writing $x = p + \veps$ with $\veps \in [0,1-p]$, this inequality and inequality~\eqref{ineq:K} imply that
\begin{align*}
	- \log \Pr(X \ge p + \veps) \
	&\ge \ m K(p,p+\veps) + \log(\veps/\sigma_{a,b}) + \log \sqrt{\pi/2} \\
	&\ge \ \frac{m \veps^2}{2(p + 2\veps/3)(1 - p - 2\veps/3)}
		+ \log(\veps/\sigma_{a,b}) + \log \sqrt{\pi/2} .
\end{align*}
Hence, it suffices to show that for any $d > 0$ and sufficiently large $c(d) > 0$,
\begin{equation}
\label{ineq:Bernstein.remainder}
	\log(\veps/\sigma_{a,b}) + \log \sqrt{\pi/2} \
	\ge \ \frac{d \veps^2}{2(p + 2\veps/3)(1 - p - 2\veps/3)}
\end{equation}
whenever $c(d) \sigma_{a,b} \le \veps \le 1 - p$. In fact, the right-hand side of \eqref{ineq:Bernstein.remainder} is increasing in $\veps$ and equals
\[
	\frac{d (1 - p)^2}{2[p + 2(1 - p)/3] [1 - p - 2(1-p)/3]} \
	= \ \frac{9d (1 - p)}{2(2+p)} \
	\le \ 9d/4
\]
if $\veps = 1-p$, while the left-hand side of \eqref{ineq:Bernstein.remainder} is at least $\log c(d) + \log\sqrt{\pi/2}$ if $\veps \ge c(d) \sigma_{a,b}$. This shows that the assertion of the corollary is true with $c(d)$ at most $\sqrt{2/\pi} \exp(9d/4)$.
\end{proof}

\begin{proof}[\bf Proof of Lemma~\ref{lem:R(sigma.delta)}]
At first we analyze $\tilde{r}(a)$. We use Binet's formula for the remainder $r(y)$ in \eqref{eq:Stirling},
\[
	r(y) \ = \ \int_0^\infty e^{-yt} w(t) \, \d t
\]
with a certain function $w$ satisfying $12^{-1} e^{-t/12} < w(t) < 12^{-1}$, see D{\"u}mbgen et al.\ (2021, Lemma~10). Consequently,
\[
	2 r(a) - r(2a) \ = \ \int_0^\infty (2e^{-at} - e^{-2at}) w(t) \, \d t ,
\]
and since $2e^{-at} - e^{-2at} = e^{-at} (2 - e^{-at}) > 0$, we conclude that
\[
	2 r(a) - r(2a) \
	\begin{cases}
		\displaystyle
		< \ \frac{1}{12} \int_0^\infty (2e^{-at} - e^{-2at}) \, \d t
			\ = \ \frac{1}{8a} , \\[2ex]
		\displaystyle
		> \ \frac{1}{12} \int_0^\infty (2 e^{-(a+1/12)t} - e^{-(2a+1/12)t}) \, \d t
			\ = \ \frac{a+1/36}{8 (a + 1/12)(a + 1/24)} .
	\end{cases}
\]
In particular, as $a \to \infty$,
\begin{equation}
\label{eq:rtilde.asymptotic}
	\tilde{r}(a) \ = \ - \frac{1}{8a} + O(a^{-2}) .
\end{equation}
Moreover,
\begin{equation}
\label{ineq:rtilde.derivative}
	\frac{d}{da} \tilde{r}(a) \
	= \ 2 \int_0^\infty t (e^{-at} - e^{-2at}) w(t) \, \d t \
	< \ \frac{1}{6} \int_0^\infty t (e^{-at} - e^{-2at}) \, \d t \
	= \ \frac{1}{8a^2} .
\end{equation}

Next we verify that $\log R(\sigma(a))$ is strictly decreasing in $a > 1$. It follows from representation \eqref{eq:R.sigma.2} that
\begin{align}
\label{eq:R.sigma.3}
	\log R(\sigma(a)) \
	= \ &\tilde{r}(a) + \log(a)/2 \\
\nonumber
	&+ \ \delta - (a - 1/2) \log(a+\delta)
		+ 1 + (a-1) \log(a-1) .
\end{align}
According to \eqref{ineq:rtilde.derivative}, the derivative of this with respect to $a$ is not greater than
\[
	\frac{1}{8a^2} + \frac{1}{2a} - \frac{a-1/2}{a+\delta} - \log(a + \delta)
		+ \log(a-1) + 1 .
\]
For fixed $a > 1$, the derivative of this bound with respect to $\delta \ge 1$ equals $- (\delta + 1/2) /(a + \delta)^2$, so it is maximal for $\delta = -1/2$. This leads to
\begin{align*}
	\frac{d}{da} \log R(\sigma(a)) \
	&\le \ \frac{1}{8a^2} + \frac{1}{2a} + \log \Bigl( \frac{a-1}{a-1/2} \Bigr) \
		 = \ \frac{1}{8a^2} + \frac{1}{2a} + \log \Bigl( 1 - \frac{1}{2(a - 1/2)} \Bigr) \\
	&< \ \frac{1}{8a^2} + \frac{1}{2a} + \log \Bigl( 1 - \frac{1}{2a} \Bigr) \
		= \ - \sum_{k\ge 3} (2a)^{-k}/k \ < \ 0 .
\end{align*}

It remains to prove the expansion of $\log R(\sigma(a)))$ as $a \to \infty$. To this end, we rewrite \eqref{eq:R.sigma.3} as
\[
	\log R(\sigma(a)) \
	= \ \tilde{r}(a) + \delta - (a-1/2) \log(1 + \delta/a) + (a-1) \log(1 - 1/a) .
\]
Since $\log(1 + y) = y + O(y^2) = y - y^2/2 + O(y^3)$ as $y \to 0$,
\begin{align*}
	\delta - (a-1/2) \log(1 + \delta/a) \
	&= \ \delta - \frac{(a-1/2)\delta}{a} + \frac{(a-1/2) \delta^2}{2a^2} + O(a^{-2}) \\
	&= \ \frac{\delta(\delta+1)}{2a} + O(a^{-2}) , \\
	1 + (a-1) \log(1 - 1/a) \
	&= \ 1 - \frac{a-1}{a} - \frac{a-1}{2a^2} + O(a^{-2}) \\
	&= \ \frac{1}{2a} + O(a^{-2}) .
\end{align*}
Combining these expansions with \eqref{eq:rtilde.asymptotic} leads to the desired expansion of $\log R(\sigma(a))$.
\end{proof}

\paragraph{Acknowledgements.}
We are grateful to Maciej Skorski and two anonymous referees for constructive comments on earlier versions. This work was supported by the Swiss National Science Foundation.


\end{document}